\documentclass[amsmath,amssymb,11pt,tightenlines]{revtex4-2}

\usepackage{algorithm}% http://ctan.org/pkg/algorithm
\usepackage{algpseudocode}% http://ctan.org/pkg/algorithmicx
\usepackage{enumitem}
\usepackage{graphicx}

\newtheorem{theorem}{Theorem}[section]

\newcommand{\mathd}{\mathrm{d}}

\newcommand{\mathe}{\mathrm{e}}
\newcommand{\imag}{\mathrm{i}}

\newcommand{\Tinlet}{T_{\mathrm{inlet}}}
\newcommand{\Toutlet}{T_{\mathrm{outlet}}}
\newcommand{\Tinit}{T_{\mathrm{init}}}
\newcommand{\Tref}{T_{\mathrm{ref}}}

\newcommand{\mypi}{\mathcal{I}}

\newcommand{\lagrange}{L}
\newcommand{\length}{\ell}

\graphicspath{{figs/}}

\usepackage{tikz}

\usetikzlibrary{patterns}
\usetikzlibrary{decorations.pathmorphing,shapes}

\begin{document}

\title{A Mathematical Framework for Spatio-Temporal Control in Industrial Drying}

\author{Lennon \'O N\'araigh}
\email{Corresponding author.  Email: onaraigh@maths.ucd.ie}
\address{School of Mathematics and Statistics, University College Dublin, Belfield, Dublin 4, Ireland}
\date{today}

\begin{abstract}
We introduce two models of industrial drying - a simplified one-equation model, and a detailed three-equation model.  The purpose of the simplified model is rigorous validation of numerical methods for PDE-constrained optimal control.  The purpose of the detailed model is to be able to predict and control the behaviour of an industrial disk drier.  For both models,  we introduce a fully validated numerical method to compute the optimal source term to maintain the outlet temperature as close as possible to the set-point temperature.  By performing simulations using realistic parameters for industrial driers, we illustrate potential applications of the method.
%For the detailed model, the optimal control theory establishes a target for the lowest possible root-mean-square temperature fluctuation achievable at the outlet.  
%This then enables us to assess the performance of various practical control techniques, and to quantify how close they are to to the optimal control.
\end{abstract}


\maketitle

\section{Introduction}
\label{sec:intro}

In this article, we consider the mathematical modelling of  industrial driers.  Such machines  are used in the food-processing industry to remove moisture from raw materials; they require large quantities of energy to raise the temperature of the moisture and to evaporate the moisture in the raw materials.  The  machine we have in mind for modelling is the disk dryer, however, the mathematical framework we develop can be applied to a wide variety of driers and heaters.  Such machines do not always operate in steady state, as the stream of input raw materials may vary in temperature and composition.  Therefore, in order for the final moisture level and temperature at the outlet of the machine to have a fixed value, a heat source must be present to deliver a controlled amount of thermal energy to the system, to correct for the inlet fluctuations and hence, to `steer' the system to a target state at the outlet.  Mathematical models for industrial drying can give important insights into the correct way to control such systems.  Hence, the aim of this work is to develop such models, with a view to formulating a control theory to maintain the outlet temperature at a specific set-point.  

\subsection{Literature Review}

The simplest possible models in the drying literature are \textit{lumped-parameter models}.  Here, either the spatial variations or the temporal variations are averaged out, leading to a set of simple ordinary differential equations based on mass- and energy-balances\cite{mujumdar2006handbook}.    Examples include Reference~\cite{hernandez2014modeling} for a disk dryer, and Reference~\cite{luz2010dynamic} for a rotary dryer.  In these cases, the spatial degrees of freedom are averaged out, leading to a system of ordinary differential equations in the time variable.  The authors then use model predictive control to maintain the outlet variables at a set point.

And yet, often, the temperature and moisture profiles along the length of the machine are important -- for instance, it is desirable for the outlet temperature to be at a set point but also, for the temperature inside the drier to stay within certain bounds, for energy efficiency but also, for food hygiene.  For this reason, Reference~\cite{iguaz2003mathematical} looks at a set of coupled ordinary differential equations for a rotary drier, with each set of equations representing a separate segment of the rotary drier.  The authors look at how the temperature profile along the drier depends on inlet conditions and flow rate -- effectively the reverse of the previous examples where now the spatial variable (i.e. the distance $x$ from the inlet, as opposed to time) is the key variable of interest.  The authors also look at how the system responds to a step change in the dryer speed.

In the present work, both the temporal and spatial variables are important.  Hence, we look at a partial-differential-equation (PDE) model for industrial driers.  In this way, the hyperbolic nature of such systems becomes apparent -- information from the inlet is convected to the outlet.  Hence, controlling the moisture level and the temperature at the outlet requires the propagation of information from the inlet, as well as manipulation of the heat source inside the drier, to `steer' the system towards a desired outlet state.

This introduces the idea of optimal control theory.  Optimal control theory is well understood in the context of ordinary differential equations, where it is required to minimize a cost functional based on state variables, subject to the constraint that the state variables satisfy a given mathematical model.  
A time-dependent forcing term (the control variable) provides the necessary additional degrees of freedom so that the minimum can be realised -- the control variable which realises the minimum of the cost functional is called the optimal control.  

For ordinary differential equations, the optimal control can be computed using Pontryagin's Maximum Principle\cite{evans1983introduction} -- the idea here is to introduce a constrained penalty function comprising the basic objective function to be minimized, as well as a set of constraints (parametrized by Lagrange multipliers) enforcing the condition that the state variables satisfy the mathematical model.  By taking variations of the constrained functional, one obtains a set of forward equations (the  mathematical model itself), as well as a set of adjoint equations (involving the Lagrange multipliers), the solution of which gives the optimal control. 

In the present case, we apply optimal control theory to partial differential equations, where a similar technique applies.
The reader is referred to Reference~\cite{hinze2008optimization} for details.  The same technique has been used in a PDE model for the optimal control of chemical reactors in microchannels\cite{blauth2021optimal}, a PDE model for the optimal draw speed in the production of glass (`tube drawing')\cite{butt2012optimal}, and PDE model for the optimal air temperature / air speed for melt spinning, a key manufacturing technique in the production of synthetic textile fibers\cite{gotz2010optimal}.  The technique has also been applied to Burgers equation\cite{ou2011unsteady}, a canonical non-steady problem and similar in its mathematical form to the models considered in the present work.

\subsection{Plan of the paper}

References~\cite{blauth2021optimal} and~\cite{butt2012optimal} use generalized Newton methods for the minimization of the constrained cost functional and the computation of optimal controls, an approach we also use in this paper.  However, our own first steps are more cautious, and we start with a simpified model of industrial heating (Section~\ref{sec:simple}).  Crucially, the simplified model admits exact analytical solutions -- including for the  optimal control.  
In Section~\ref{sec:simple1}, we introduce a numerical optimal-control theory for the simplified model, and we use the available  analytical solutions to validate the numerical method rigorously.

Next, we introduce a more detailed model of an industrial drier, which convects the product from the inlet to the outlet (Section~\ref{sec:drier_theory}).  The product is made up of two streams: dry solid and liquid.  The model allows for the evaporation of the liquid and also involves the temperature of the streams -- hence, a three-equation model.  
We  look at the equilibrium solutions of the model, corresponding to steady inlet conditions, as well as the linear stability of such equilibrium solutions (Section~\ref{sec:linear}).  We investigate small-amplitude perturbations away from equilibrium, and demonstrate analytically that these perturbations can be suppressed so that the outlet temperature remains at a set point.  We follow this analysis with steepest-descent calculations to determine numerically the optimal heat source to maintain the outlet temperature at a fixed point.  We do this for both the case of small-amplitude inlet perturbations (linear theory, Section~\ref{sec:linear}), as well as large-amplitude perturbations, which involves a non-linear theory (Section~\ref{sec:nonlinear}).    Discussion and concluding remarks are presented in Section~\ref{sec:conc}.

\subsection{Scope of the work}

The three-equation model has been formulated with disk driers in mind.  The simplified mathematical model is more generic.  Its main purpose is a test-bed for the optimal-control theory and a means of performing rigorous validation tests on our numerical methods.  Sections~\ref{sec:simple}--\ref{sec:simple1} on the simplified mathematical model also help us to formulate the optimal-control theory in a clear and transparent manner.  Then, it is possible to describe the optimal control theory for the three-equation model by simple analogy.  The main variable targeted for control in the present work is the outlet temperature.  In practical settings, it is important to maintain this variable near a set-point for food-hygiene and quality-control reasons.  Clearly, in an industrial drier, the outlet moisture is a key variable whose value needs to be controlled.  In the present work, this variable is maintained at the equilibrium value by the choice of design parameters -- in particular, the drier length and the drier residence time.  Introducing a second manipulated variable (in addition to the heat source) with a view to gaining more control over the outlet moisture level will be the subject for future work, for which the mathematical framework introduced herein can be used.

\section{Simplified Mathematical Model}
\label{sec:simple}

We look at a simplified mathematical model of drying.  This model involves only the product temperature, and reads as follows:
\begin{subequations}
\begin{equation}
\frac{\partial T}{\partial t}+u_0\frac{\partial T}{\partial x}=k\left[q(t)-T\right],\qquad x\in (0,\length].
\label{eq:toy_global_pde}
\end{equation}
Here, $u_0>0$ is a constant velocity (dimensions: $[\text{Length}]/[\text{Time}]$), and $k>0$ is a rate (dimensions: $1/[\text{Time}]$).  The quantity $q(t)$ (with dimensions of temperature) represents the temperature of the surroundings, this may be manipulated to bring $T(x,t)$ to a desired value at the outlet $x=\length$ (i.e. $\Toutlet(t)=T(\length,t)$).
%
%Manipulating $q(t)$ to control $\Toutlet(t)$ is of practical concern as the product temperature at the inlet varies in time as different materials are processed through the drier.  
Hence, the  temperature at the inlet varies in time according to:
\begin{equation}
T(0,t)=\Tinlet(t),
\label{eq:toy_global_inc}
\end{equation}%
where $\Tinlet(t)$ is a time-dependent inlet condition, with $t\geq 0$.  Furthermore, at $t=0$, the temperature is determined by an initial condition:
\begin{equation}
T(x,0)=\Tinit(x),
\label{eq:toy_global_ic}
\end{equation}%
with $x\geq 0$.  The solution is continuous if $\Tinit(0)=\Tinlet(0)$.  
\label{eq:toy_global}%
\end{subequations}%

\subsection{Analytical and Numerical Solutions}

Equation~\eqref{eq:toy_global} admits an analytical solution, which may be validated by direct computation:
\begin{subequations}
\begin{equation}
T(x,t)=\begin{cases} \Tinit(x-u_0 t)\mathe^{-kt}+\mypi(t),&u_0 t<x,\\
                    \left[\Tinlet\left(t-\frac{x}{u_0}\right)-\mypi\left(t-\frac{x}{u_0}\right)\right]\mathe^{-kx/u_0}+\mypi(t),&u_0 t>x,
										\end{cases}
\end{equation}
where 
%$\mypi(t)$ is the particular integral for the spatially homogenous problem $\mathd \mypi/\mathd t=k[q(t)-\mypi]$, with $\mypi(0)=0$:
%
\begin{equation}
\mypi(t)=\mathe^{-kt}\int_0^t \mathe^{kt'}[k\,q(t')]\,\mathd t'.
\end{equation}%
\label{eq:asln}%
\end{subequations}%
Equation~\eqref{eq:asln} clearly highlights the hyperbolic nature of the  model in Equation~\eqref{eq:toy_global}: at early times ($u_0 t<x$), information from the initial condition dominates, and is carried downstream by the constant velocity $u_0$.  At later times ($u_0 t>x$), information from the inlet dominates, and is carried downstream by the flow.  
%As we have chosen initial and inlet conditions carefully, the solution is continuous along the characteristic curve $x=u_0 t$: $T(x=u_0t ,t)=\Tinit(0)\mathe^{-kt}=\Tinlet(0)\mathe^{-k(u_0 t)/u_0}$.

%\subsection{Numerical Method}

 Equation~\eqref{eq:toy_global} can also be described numerically.  For this purpose,  the equation is discretized in space and time using a second-order upwind scheme:
\begin{equation}
\frac{T_i^{n+1}-T_i}{\Delta t}+u_0\left(\frac{3T_i^n-4T_{i-1}^n+T_{i-2}^n}{2\Delta x}\right)=
k\left[q(t^n)-T_i^n\right].
\label{eq:numsln}%
\end{equation}%
This is an explicit scheme, as such it is constrained by a CFL condition, $u_0\Delta t<\Delta x$.  The numerical approach is useful in caess where analytical solutions are not available, for instance, in the three-equation model.

%Here, we use standard notation, where time is discretized as $t^n=n\Delta t$, where $n=0,1,2,\cdots$ is an integer labelling the current timestep and $\Delta t$ is the stepsize.  Similarly, for the spatial discretization, we use $x_i=i\Delta x$, where $i=0,1,\cdots,N-1$, and $\Delta x=L/(N-1)$ is the grid size in the spatial domain.  Hence, $T_i^n=T(x_i,t^n)$.

\subsection{Optimal Control -- Special Analytical Solution}

For the Simple Mathematical Model, an analytical solution for the optimal control $q(t)$ can be constructed, which maintains the outlet temperature $T(\length,t)$ at a set-point $T_*$.  To see this,  we simply set $T(\length,t)=T_*$
in the analytical solution in Equation~\eqref{eq:asln}.  This gives:
\begin{subequations}
\begin{equation}
T_*=\Tinit(\length-u_0 t)\mathe^{-kt}+\mypi(t),\qquad u_0 t<\length,
\label{eq:oc_allb}%
\end{equation}
and
\begin{equation}
T_*=\left[\Tinlet(t-\length/u_0)-\mypi(t-\length/u_0)\right]\mathe^{-k\length/u_0}+\mypi(t),\qquad u_0 t > \length.
\label{eq:oc_allc}%
\end{equation}%
\label{eq:oc_all}%
\end{subequations}%
Multiply Equation~\eqref{eq:oc_allb} by $\mathe^{kt}$ and differentiate with respect to time to obtain:
\begin{subequations}
\begin{equation}
T_*+\frac{u_0}{k}\Tinit'(\length-u_0 t)\mathe^{-kt}+q(0)\mathe^{-kt}=q(t),\qquad t<\length/u_0.
\label{eq:oc_all_sln_a}
\end{equation}
A similar treatment of Equation~\eqref{eq:oc_allc} gives:
\begin{multline}
T_*-\Tinlet(t-\length/u_0)\mathe^{-k\length /u_0}-\frac{1}{k}\Tinlet'(t-\length/u_0)\mathe^{-k\length/u_0}\\
+q(t-\length/u_0)\mathe^{-k\length/u_0}=q(t),\qquad t>\length/u_0.
\label{eq:oc_all_sln_b}
\end{multline}%
\label{eq:oc_all_sln}%
\end{subequations}%
Equations~\eqref{eq:oc_all_sln} contain an unknown quantity $q(0)$.  However, in the case where the initial temperature, the mean temperature $T_0$, and the set-point temperature $T_*$ are all the same ($\Tinit=T_0=T_*$), $q(0)$ can be taken as zero.  Then, Equation~\eqref{eq:oc_all_sln_a} has a closed-form solution, and Equation~\eqref{eq:oc_all_sln_b} can be solved recursively.  Notice that the optimal control has a possible jump discontinuity at the characteristic time $ t_0=\length/u_0$, with
\begin{eqnarray*}
q(t_0^-)&=&T_*+\frac{u_0}{k}\Tinit'(0)\mathe^{-kt_0},\\
q(t_0^+)&=&T_*-\Tinlet(0)\mathe^{-kt_0}-\frac{1}{k}\Tinlet'(0)\mathe^{-kt_0}.
\end{eqnarray*}
%	

%The optimal control $q(t)$ is plotted in Figure~\ref{fig:qt1} (parameter values the same as before).
%\begin{figure}[htb]
	%\centering
		%\includegraphics[width=0.5\textwidth]{qt1}
		%\caption{The optimal control $q(t)$ (Equation~\eqref{eq:oc_all_sln})}
	%\label{fig:qt1}
%\end{figure}
%%
%%
%The results of applying this  control to the discretized mathematical model~\eqref{eq:numsln} are shown in Figure~\ref{fig:Ttoy_xt_oc} (parameter values the same as before) -- the outlet set-point condition $T(\length,t)=100^\circ\mathrm{C}$ is attained.  
%%	
%\begin{figure}[htb]
	%\centering
		%\includegraphics[width=0.5\textwidth]{Ttoy_spacetime_oc}
		%\caption{Plot showing the behavior of the numerical solution in space and time, with the optimal control $q(t)$ (Equation~\eqref{eq:oc_all_sln}) applied.}
	%\label{fig:Ttoy_xt_oc}
%\end{figure}
%%
%%

\subsection{Optimal Control -- limit of large $k\length/u_0$}

Assume a bounded control, $|q(t)|\leq M$, for all $t\in[0,\infty)$, and a bounded inlet condition $|\Tinlet(t)|\leq M'$ for all $t\in[0,\infty)$.  Then, $\mypi(t)$ is bounded by $\mypi(t)\leq kM(1-\mathe^{-kt})$.  Hence, in the limit of large $k\length/u_0$, the analytical solution~\eqref{eq:asln} simplifies:
\begin{equation}
T(\length,t)=\begin{cases} \Tinit(\length-u_0 t)\mathe^{-kt}+\mypi(t),& t<\length/u_0,\\
                    \mypi(t),& t>\length/u_0,
										\end{cases}
\end{equation}
Hence, the influence of the inlet at late times $t>\length/u_0$ is `washed away' by the exponential factor $\mathe^{-k\length/u_0}$.  Thus, at late times, the outlet temperature can be manipulated by choice of $\mypi(t)$.  Simply taking $q(t)=T_*$ gives $\mypi(t)=T_*(1-\mathe^{-kt})$ and hence,
\begin{equation}
T(\length,t)=T_*(1-\mathe^{-kt}),\qquad t>\length/u_0,
\end{equation}
Thus, the outlet condition at late times is controlled simply by taking $q(t)=T_*$, provided $k\length/u_0\gg 1$.

\section{General Optimal Control Theory for the Simplified Mathematical Model}
\label{sec:simple1}

To extend our method of control to cases where analytical solutions are not possible, we formulate a more general theory of optimal control based on a steepest-descent technique.
Hence, we introduce the penalty function
\begin{equation}
J(T,q)=\tfrac{1}{2}\int_0^\tau \left[T(\length,t)-T_*\right]^2\mathd t.
\end{equation}
Here, $\tau$ is the final simulation time.
The aim of this section is to develop a general framework for computing the optimal control $q(t)$ at each point in time  $t\in [0,\tau]$, such that $J$ is minimized.

For these purposes, we seek to minimize the functional $J$ subject to the constraint that $T(x,t)$ satisfies the mathematical model~\eqref{eq:toy_global}.  Hence,   we introduce a
 a constrained functional:
\begin{equation}
\lagrange(T,q,\psi)=\tfrac{1}{2}\int_0^\tau \left[T(\length,t)-T_*\right]^2\mathd t+
\int_0^\tau\mathd t\int_0^\length \mathd x
\left[\frac{\partial T}{\partial t}+u_0\frac{\partial T}{\partial x}
-k\left(q(t)-T\right)\right]\psi.
\label{eq:lagrange1}
\end{equation}
Here, $\psi$ is a Lagrange multiplier.  The physical interpretation fo $\psi$ is that it enforces the constraint that $T(x,t)$  satisfies the model equations (in this context, Equation~\eqref{eq:toy_global_pde}).   We follow standard practice in developing an equation for $\psi$.  However, we note that in some applications, the vanishing of the Lagrange multiplier engenders the so-called Lagrange catastrophe, which may cause the constraints to be violated~\cite{he2020lagrange}.  In this case, the semi-inverse method can be used to develop the appropriate functional to which a variational principle can be applied~\cite{wang2020variational,he2020lagrange}.

We next compute $\delta \lagrange$, which is the variation in $\lagrange$ with respect to variations in $T$, and variations in $q$ (denoted by $\delta q$).  Note that $\delta q$ depends only on time, not on space.  We compute $\delta \lagrange$ in a straightforward fashion and obtain:
\begin{multline}
\delta \lagrange=\int_0^\tau\left[T(\length,t)-T_*\right]\delta T(\length,t)\mathd t
-\int_0^\tau \mathd t\int_0^\length k\psi\delta q\,\mathd x\\
+\int_0^\tau \int_0^\length \left[-\frac{\partial\psi}{\partial t}-u_0\frac{\partial\psi}{\partial x}+k\psi\right]\delta T\,\mathd x+
\int_0^\length\left[\psi\delta T\right]_{t=\tau}\mathd x+\int_0^\tau \left[u_0\psi\delta T\right]_{x=\length}\mathd t.
\end{multline}
Here, in computing the boundary terms (which arise from using integration by parts), we have assumed that $\delta T=0$ at the inlet and also, $\delta T=0$ at $t=0$.  This makes sense, since inlet/initial conditions are fixed, so there is no variation in the solution at these points.  Based on this expression for $\delta \lagrange$, we identify the adjoint problem:
\begin{subequations}
\begin{equation}
-\frac{\partial \psi}{\partial t}-u_0\frac{\partial \psi}{\partial x}+k\psi=0,\qquad x\in [0,\length).
\end{equation}
We impose the outlet condition 
\begin{equation}
\psi(\length,t)=-(1/u_0)[T(\length,t)-T_*],
\end{equation}
and we further impose the terminal condition 
\begin{equation}
\psi(x,\tau)=0.
\end{equation}%
\label{eq:adjoint}%
\end{subequations}%
With the adjoint problem defined in this way, the variation in $\delta \lagrange$ simplifies dramatically:
\begin{equation}
\delta \lagrange=-\int_0^T \mathd t \,\delta q\int_0^\length k\psi(x,t)\mathd x.
\end{equation}
These calculations now provide a numerical method for the computation of the optimal control $q(t)$:
\begin{algorithm}[H]
	\label{algo:newton}
  \caption{Steepest-Descent Algorithm for the Simplified Mathematical Model}
	\begin{algorithmic}[1]
  \State Initialize the control using some initial guess $q(t)=q_0(t)$.
	\State Solve Equation~\eqref{eq:toy_global} forwards from $t=0$ to the final time $\tau$, compute $T(\length,t)$ for each $t\in[0,\tau]$.
	\State Solve the adjoint problem~\eqref{eq:adjoint} backwards, from $t=\tau$ to $t=0$.
	\State Update the control $q(t)\leftarrow q(t)+\lambda\int_0^\length\left(k\psi\right)\,\mathd x$, where $\lambda$ is a small positive parameter.
	\State Repeat steps 2--4 until $L$ is sufficiently small.
	\end{algorithmic}
\end{algorithm}
%
%\begin{enumerate}[noitemsep]
%\item Initialize the control using some initial guess $q(t)=q_0(t)$.
%\item Solve Equation~\eqref{eq:toy_global} to the final time $\tau$, compute $T(\length,t)$ for each $t\in[0,\tau]$.
%\item Solve the adjoint problem~\eqref{eq:adjoint}.
%\item Update the control $q(t)\rightarrow q(t)+\lambda\int_0^\length\left(k\psi\right)\,\mathd x$, where $\lambda$ is a small positive parameter.
%\item Repeat steps 2--4.
%\end{enumerate}
%
\noindent Operations 2--4 correspond to one iteration in a cycle.  Label the iterations by $n$, with $n=0,1,2,\cdots$.  By carrying out repeated iterations, we obtain:
\begin{equation}
\delta \lagrange=\lagrange^{n+1}-\lagrange^n=-\lambda^n\int_0^\tau \left[\int_0^\length  \left(k\psi^n\right)\,\mathd x\right]^2\mathd t,
\end{equation}
for sufficiently small $\lambda^n$, and correspondingly,
\begin{equation}
q^{n+1}(t)=q^n(t)+\lambda^n \int_0^\ell (k\psi^n)\mathd x.
\label{eq:VIM}
\end{equation}
Hence, $\delta \lagrange<0$ and the iterative process will converge, such that $\delta \lagrange\rightarrow 0$ as $n\rightarrow \infty$, such that $\lagrange$ 
becomes stationary, and the optimal control for $q(t)$ is obtained.  Indeed, this algorithm is essentially a  steepest-descent method to compute the minimum of the functional $\lagrange$, where the search direction is given by:
\begin{equation}
-\frac{\delta \lagrange}{\delta q}(t)=\int_0^\length k\psi(x,t)\mathd x.
\label{eq:optimality}
\end{equation}
With the steepest-descent method identified in this manner, it is possible to use standard optimization techniques to compute the optimal value of $\lambda^n$ at each iteration cycle, this is addressed later in the section.

Finally, we notice that Equation~\eqref{eq:VIM} is reminiscent of the Variational Iteration Method for solving non-linear ordinary differential equations~\cite{he2007variational}.  The Variational Iteration Method has also been used to solve problems in optimal control theory for ordinary differential equations~\cite{berkani2012optimal} and also, partial differential equations~\cite{akkouche2014optimal}.  Hence, the Variational Iteration Method may be applicable to the present optimal-control problem as dwell, and would be complementary to the numerical method developed herein.

\subsection{Mathematical Setting}

We place the foregoing calculations in a more mathematical context as follows.  We let $T(x,t)\in Y$ and $q(t)\in U$.  We let $e(T,q)=\partial_t T+u\partial_x T-k(q-T)\in Z$.  As $e$ is linear in $T$ and $q$, this can also be identified with a linear problem,
\begin{equation}
e(T,q)=AT+Bq,\qquad A\in \mathcal{L}(Y,Z).
\end{equation}
Here, $Y$ and $Z$ are appropriate Banach spaces, $\mathcal{L}(Y,Z)$ is the space of linear operators from $X$ to $Y$, and $U$ is an appropriate Hilbert space.  Hence, the weak form of the simplified mathematical model reads:
\begin{equation}
e(T,q)=0\text{ in }Z,
\end{equation}
where $e: Y\times U\rightarrow Z$ is defined with respect to its action on test functions:
\begin{multline}
\langle \psi,e(T,q)\rangle_{Z^*,Z}=
\int_0^\tau \mathd x\int_0^\length (-T\partial_t \psi)\,\mathd x\\
+\int_0^\tau \mathd x\int_0^\length u(-T\partial_x \psi)\,\mathd x
-\int_0^\tau \mathd x\int_0^\length k(T-q)\psi\,\mathd x\\
+\int_0^\tau u\left[T(\length,t)\psi(\length,t)-T(0,t)\psi(0,t)\right]\mathd t
+\int_0^\length \left[T(x,\tau)\psi(x,\tau)-T(x,0)\psi(x,0)\right]\mathd x,
\end{multline}
for all test functions $\psi\in Z^*$.
The Lagrangian of the system is $L:Y\times U\times Z^*\rightarrow \mathbb{R}$, where
\begin{equation}
L(T,q,\psi)=J(T,u)+\langle \psi,e(T,q)\rangle_{Z^* Z}
\end{equation}
which we recognize as Equation~\eqref{eq:lagrange1}.  The first-order optimality conditions in the direction $(\delta T,\delta\psi,\delta u)$ are:
\begin{subequations}
\begin{eqnarray}
%J_T(T,u)\delta T+\langle \lambda,e_T(T,q)\delta T\rangle_{Z^*,Z}&=&0,\label{eq:abs1}\\
%\langle \delta\lambda,e_T(T,q)\delta T\rangle_{Z^*,Z}&=&0,\label{eq:abs2}\\
%J_q(T,q)\delta q+\langle \lambda,e_q(T,q)\delta q\rangle_{Z^*,Z}&=&0;\label{eq:abs3}
\langle J_T(T,q),\delta T\rangle_{Y^*,Y}+\langle \psi,e_T(T,q)\delta T\rangle_{Z^*,Z}&=&0,\label{eq:abs1}\\
\langle \delta\psi,e_T(T,q)\delta T\rangle_{Z^*,Z}&=&0,\label{eq:abs2}\\
\langle J_q(T,q),\delta q\rangle_{U^*,U}+\langle \psi,e_q(T,q)\delta q\rangle_{Z^*,Z}&=&0;\label{eq:abs3}
\end{eqnarray}
\end{subequations}
these conditions must hold for all (admissible) directions $\delta T\in Y$, $\delta \psi\in Z$, and $\delta q\in U$.  We recognize Equation~\eqref{eq:abs1} as the adjoint problem~\eqref{eq:adjoint}, Equation~\eqref{eq:abs2} as the weak formulation of the simplified mathematical model (Equation~\eqref{eq:toy_global_pde}), and Equation~\eqref{eq:abs3} as the optimality condition $\delta J/\delta q=0$ (\textit{cf.} Equation~\eqref{eq:optimality}).

Furthermore, we look at solutions to $e(T,q)=0$, such that $T$ depends implicitly on the forcing term, $T=T(q)$.  Hence, we look at the reduced cost function,
\begin{equation}
\widehat{J}(q)=J(T(q),q)=J(T(q),q)+\langle \psi,e(T(q),q)\rangle_{Z^*,Z}=L(T(q),q,\psi).  
\end{equation}
Differentiating this expression in a direction $\delta q\in U$, we obtain:
\begin{equation}
\langle \widehat{J}'(q),\delta q\rangle_{U^*,U}=\langle L_T(T(q),q,\psi),\frac{\partial T}{\partial q}\delta q\rangle_{Y^* Y}+\langle L_q(T(q),q,\psi),\delta q\rangle_{U^*,U}.
\end{equation}
We now choose $\psi=\psi(q)$ such that:
\begin{equation}
L_T(T(q),q,\psi)=0
\end{equation}
(this is nothing but the adjoint equation~\eqref{eq:abs1}).  
In fact,
\begin{eqnarray*}
\langle L_T(T,q,\psi),d\rangle_{Y^*,Y}&=&\langle J_T(T,q),d\rangle_{Y^*,Y}+\langle \psi,e_T(T,q)d\rangle_{Z^*,Z},\\
                                         &=&\langle J_T(T,q)+e_T(T,q)^*\psi,d\rangle_{Z^*,Z}.
\end{eqnarray*}
Furthermore, we now have:
\begin{eqnarray*}
\langle \widehat{J}'(q),\delta q\rangle_{U^*,U}&=&\langle L_q(T(q),q,\psi(q)),\delta q\rangle_{U^*,U},\\
&=&\langle J_q(T(q),q),\delta q\rangle_{U^*,U}+\langle \psi(q),e_q(T(q),q),\delta q\rangle_{U^*,U},\\
&=&\langle J_q(T(q),q),\delta q\rangle_{U^*,U}+\langle e_q(T(q),q)^*\psi(q),\delta q\rangle_{U^*,U},
\end{eqnarray*}
which is true for all admissible perturbations $\delta q\in U$.  As $U$ is a Hilbert space, the Riesz representation theorem applies, and we can identify the vector $\widehat{J}'(q)\in U^*=U$:
\begin{equation}
\widehat{J'}(q)=J_q(T(q),q)+e_q(T(q),q)^*\psi(q),
\end{equation}
which is precisely Equation~\eqref{eq:optimality}.

In Section~\ref{sec:simple}, we have already established (by direct calculation) that the  optimization problem
\begin{equation}
\min_{(T,q)\in Y\times U}J(T,q),\qquad \text{subject to } AT+Bq=0,\qquad q\in U,\qquad T\in Y,
\end{equation}
has a minimum -- specifically, that $\min_{(T,q)\in Y\times U}J(T,q)=0$.   We also look here at the more general case,
\begin{multline}
\min_{(T,q)\in Y\times U}J(T,q)+\tfrac{1}{2}\alpha \|q\|_U^2,\text{subject to } AT+Bq=0,\\ u\in U_{ad}\subset U,\qquad T\in Y_{ad}\subset Y,
\label{eq:optcontgen1}
\end{multline}
Then, the following result holds\cite{hinze2008optimization}:
%thu
\begin{theorem}  {\textbf{(Hinze et al.)}}  Suppose the following statements are true:
\begin{enumerate}[noitemsep]
\item $\alpha \geq 0$, $U_{ad}\subset U$ is convex, and in the case of $\alpha=0$, bounded;
\item $Y_{ad} \subset Y$ is convex and closed, such that $AT+Bq=0$ for at least one set of functions $T\in Y_{ad}$ and $q\in U_{ad}$;
\item $A\in \mathcal{L}(Y,Z)$ has a bounded inverse.
\end{enumerate}
Then, the problem~\eqref{eq:optcontgen1} has an optimal solution $(T_{opt},q_{opt})$.  If $\alpha>0$, then the solution is unique.
\label{thm:hinze}
\end{theorem}
For the present case, the form of $J$ suggests the choice $U=L^2([0,\tau])$.  
We limit ourselves to bounded controls, such that $q\in \overline{B(0,r)}=U_{ad}\subset U$.  We also take $Y_{ad}=Y$, 
and $A=\partial_t+u_0\partial_x+k$.   The explicit solution of $AT=kq$ given in Equation~\eqref{eq:asln} establishes the that $A$ has a bounded inverse.  Thus, the conditions of the Theorem are satisfied, and a unique optimal solution exists. 

Beyond the Simple Mathematical Model, we look at a detailed three-equation model in Sections~\ref{sec:drier_theory}--\ref{sec:nonlinear}.  We look at the linearized three-equation model in  Section~\ref{sec:linear}, to which Theorem~\ref{thm:hinze} again applies.  In this case, Fourier / Laplace theory can be used to construct an explicit solution for $T(q)$, which establishes the existence of a bounded inverse operator for the problem.  We furthermore look at the nonlinear version of the model in Section~\ref{sec:nonlinear}.  This takes the form of a set of coupled semi-linear parabolic equations; the equations can be reduced to a single semi-linear equation for temperature $T$, with a continuous non-linear term.  In this case, Reference~\cite{casas1995optimal} establishes the existence of an optimal control.  Other works which look at nonlinear differential equations from the point of view of variational calculus include References~\cite{wang2009variational,he2023variational}, and~\cite{wang2020variational}.

\subsection{Barzilai--Borwein Method}

%In this work, we use the Barzili--Borwein Method for choosing the stepsize $\lambda^n$ at each iteration of the steepest-Descent method.  

In the present work, we are interested in minimizing a cost function $\min_{x\in X} J(x)$.  In practice, this amounts to solving an equation of the type
\begin{equation}
G(x)=0,
\label{eq:Gx}
\end{equation}
in this context, $G(x)$ is the gradient of the cost function $J$.  Here, $X$ is an appropriate Hilbert space, with inner product $\langle \cdot,\cdot\rangle_X$, and with norm $\|x\|_2^2=\langle x,x\rangle_X$.  The generalized Newton method is used to find the solution $\overline{x}$ of Equation~\eqref{eq:Gx}\cite{hinze2008optimization}:
\begin{algorithm}[H]
	\label{algo:newton}
  \caption{Generalized Newton Method}
	\begin{algorithmic}
  \State Choose $x_0\in X$ (sufficiently close to the solution $\overline{x}$).
	\For{$k=0,1,2,\cdots$}
	\State Choose an invertible operator $M_k\in \mathcal{L}(X,X)$,
  \State Obtain $s_k$ by solving 
		\begin{equation}
			M_k s=-G(x_k).
		\label{eq:update}
		\end{equation}
   \State Set $x_{k+1}=x_k+s_k$.
	\EndFor
	\end{algorithmic}
\end{algorithm}

There are a number of methods for choosing the operator $M_k$.  In the steepest-descent method, $M_k$ is chosen to be a diagonal matrix, $M_k=\lambda_k^{-1}\mathbb{I}$, thus, the sequence of steps $\{x_0,x_1,\cdots,x_k,\cdots\}$ follows the path of steepest descent towards the minimum of $J(x)$.   The length of the descent step $\lambda_k$ is chosen such that $F(x_{k+1})<F(x_k)$, and the method is known to possess linear convergence.

The classical Newton method involves choosing $M_k$ to be the derivative of $G$.  In this case, convergence of the iterative method is super-linear, provided the derivative of $G$ is a continuously invertible linear operator on $X$\cite{hinze2008optimization}.  However, often it is expensive to compute the derivative of $G$, in which case the derivative of $G$ may be approximated by the Secant Method:
\begin{equation}
d_k=x_k-x_{k-1},\qquad g_k=G(x_k)-G(x_{k-1}),\qquad M_k d_k=g_k,\qquad \text{(Secant Method)}.
\end{equation}

The Barzilai--Borwein method combines the ideas of the Steepest-Descent method and the Secant Method.
Hence, the matrix $M_k$ in the Secant Method is approximated as a diagonal matrix:
\begin{subequations}
\begin{equation}
M_k=\frac{1}{\alpha_k}\mathbb{I}.
\end{equation}
The equation $M_k d_k=g_k$  is then solved in the least-squares sense:
\begin{equation}
\alpha_k=\mathrm{arg}\,\min_{\alpha>0}\|\frac{1}{\alpha}d_k-g_k\|_2^2,
\end{equation}
which gives
\begin{equation}
\alpha_k=\frac{\|d_k\|_2^2}{\langle d_k,g_k\rangle}_X.
\end{equation}
Then, the update step~\eqref{eq:update} in the Barzilai--Borwein method becomes:
\begin{equation}
x_{k+1}=x_k-\alpha_k g(x_k).
\end{equation}
\end{subequations}
%\label{eq:updateBB}
%
Although the Barzilai--Borwein method does not decrease the objective function monotonically,  extensive numerical experiments show that it performs
much better than a standard Steepest-Descent method\cite{huang2022acceleration}.  In the present work, we use the Barzilai--Borwein method exclusively, and monitor the convergence of the method (and the corresponding decrease in the objective function) carefully.

\subsection{Sample Numerical Results}

\begin{figure}[htb]
	\centering
		\includegraphics[width=0.5\textwidth]{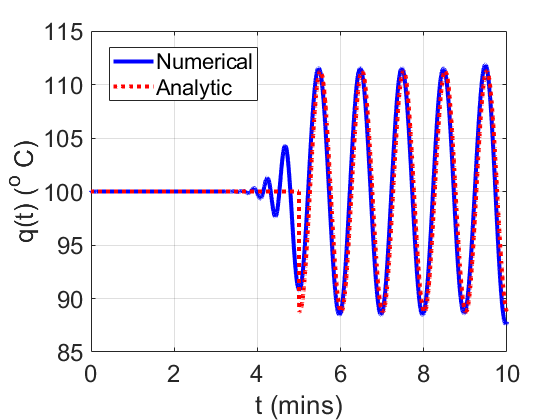}
		\caption{The optimal control $q(t)$.  Comparison between the numerical method (solid line) and the analytical solution for the optimal control (dashed line,  Equation~\eqref{eq:oc_all_sln}).}
	\label{fig:qt1_adjoint}
\end{figure}

We use the Barzilai--Borwein method to compute the optimal control $q(t)$ for the model problem~\eqref{eq:toy_global}, with $\Tinit=T_*=100^\circ\,\mathrm{C}$, time in minutes, and the other parameters given in Table~\ref{tab:params1}.  The parameter values in Table~\ref{tab:params1} are not from a real industrial drier but rather, are representative of the timescales observed in such a drier.  A detailed set of representative parameters for an industrial drier is given in Section~\ref{sec:linear}, in Table~\ref{tab:params2}.
The inlet forcing is given by:
\begin{equation}
T_{\mathrm{inlet}}(t)=\left[100+10\sin(2\pi t)\right]^\circ\mathrm{C}.
\label{eq:Tinlet}
\end{equation}
\begin{table}[htb]
			\caption{Parameter values for the simulation}
		{\begin{tabular}{|c|c|}
		\hline
		Parameter & Value\\
		\hline
			$\length$ & $5\,\mathrm{m}$\\
			$u_0$ & $1\,\mathrm{m}\cdot\mathrm{min}^{-1}$\\
			$k$   & $0.5\,\mathrm{min}^{-1}$\\
			$N$   & 200\\
			$\Delta t$ & $10^{-3}\,\mathrm{min}$\\
			\hline
		\end{tabular}}
		\label{tab:params1}
\end{table}
 A comparison between the numerical optimal control and the analytical  one  is shown in Figure~\ref{fig:qt1_adjoint}. 
Near the  characteristic time $t_0=\length/u_0$, there is some discrepancy between the two approaches.  This comes from the underlying discretization of the model partial differential equation (i.e. Equation~\eqref{eq:numsln}), which through numerical diffusion, introduces some extra smoothness into the numerical solution.  Otherwise, there is excellent agreement between the two methods, and this validates the numerical approach to computing the optimal control.

\begin{figure}[htb]
	\centering
		\includegraphics[width=0.5\textwidth]{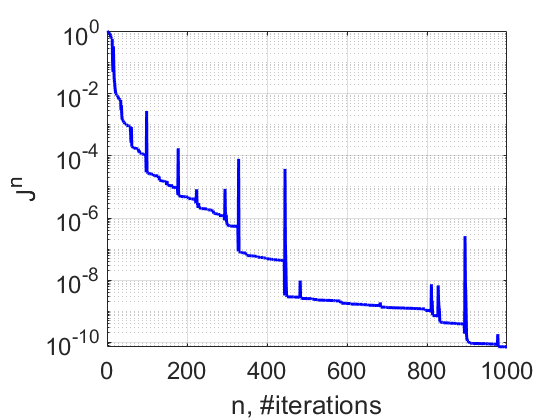}
		\caption{Plot showing the value of the cost function $J(T,q)$ at each iteration of the Barzilai--Borwein method}
	\label{fig:convergenceBB}
\end{figure}

Convergence of the Barzilai--Borwein method is shown in Figure~\ref{fig:convergenceBB}.  The convergence is non-monotone, however, each `spike' in the value of the cost function $J$ is followed by a subsequent reduction in $J$, which leads to an overall reduction in $J$ to close to a value  of $10^{-10}$ after $1000$ iterations.

\subsection{Inlet Forcing}

In the foregoing results, we have used sinusoidal inlet forcing (Equation~\eqref{eq:Tinlet}).  Often, in the food-processing industry, the inlet temperature is characterized by a stochastic process, and as such the inlet temperature consists of a mean value, plus a random perturbation.  The random perturbations can be quite complex (i.e. more complicated than a simple Brownian process), and exhibit correlations over a timescale of hours.
  We have therefore looked at using time series of temperature obtained from sensor data from real industrial driers as an inlet forcing term, both for the simplified mathematical model and the three-equation model.  The two types of inlet forcing  produce optimal controls with the same type of behaviour.  Hence, for the purpose of introducing and validating the present mathematical framework, results are presented for sinusoidal forcing throughout this article.  The applicability of the results to stochastic inlet forcing is discussed in Section~\ref{sec:conc}.

\section{Industrial Drier: Theoretical Formulation}
\label{sec:drier_theory}

\begin{figure}
	\centering
		\includegraphics[width=0.98\textwidth]{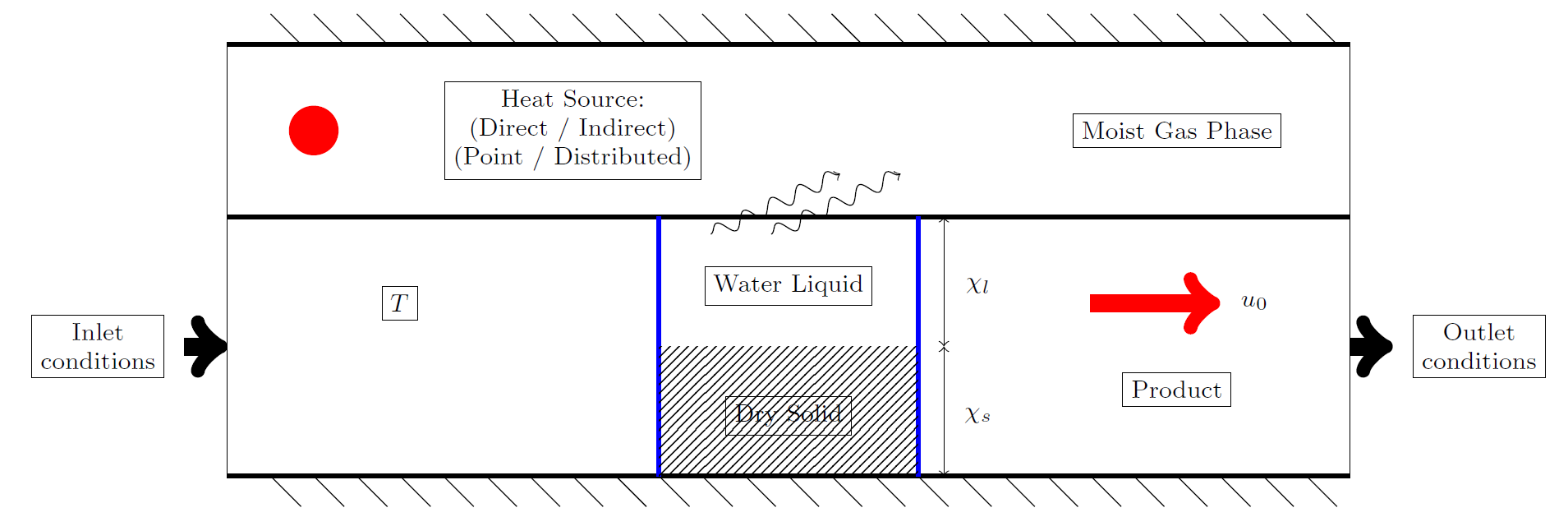}
		\caption{Schematic diagram showing the mass and energy balances for an industrial dryer}
	\label{fig:schematic1}
\end{figure}

The operating principle of an industrial disk drier is shown schematically in  Figure~\ref{fig:schematic1}.  Product is fed in to the drier at the inlet, the product consists of solids and water in the liquid phase, at a given inlet temperature $\Tinlet$.  The product is then carried through the drier at constant velocity $u_0$ and subjected to a heat source (direct contact, via heated disks in the case of the disk drier).   The product is assumed to be well mixed, such that the dry solid component and the water component are at a common temperature $T$.  The product is exposed to a moist gas phase, the heat source heats up the product and causes the water in the product to evaporate into the moist gas.  A continuous stream of air enters the drier at some location; also, the moist gas is vented from the drier.  There is therefore a continuous flow of moist gas through the drier, the flow rate is assumed to be large, such that the moist gas phase is turbulent and well mixed  and hence, the moist gas phase is treated as homogeneous and steady.
The dry solid volume fraction is denoted by $\chi_s$, the volume fraction of the water in the solid is denoted by $\chi_l$, and the volume fraction of the moist gas is $\chi_g$, with the consistency condition $\chi_s+\chi_l+\chi_g=1$.

%The system shown schematically in Figure~\ref{fig:schematic1} is suggestive of a hyperbolic system, whereby information at the inlet is propagated to the outlet.  As the inlet conditions are possibly non-steady, and as the state of product varies as it is carried along the length of the drier, the state of the components (liquid water, dry solid, and product temperature) vary in space and time. 

Assuming that the conditions of the drier are homogeneous in the plane perpendicular to the page, only the dimension along the length of the drier ($x$-direction) is important.  Hence, we identify the  three variables as determining the state of the drier at any instant: $\epsilon_s=\rho_s \chi_s$, $\epsilon_l=\rho_l \chi_l$, and $T$, the product temperature.  The internal energy $e$ per unit mass of product is a function of these variables, $e=e(\epsilon_s,\epsilon_l,T)$.  Conservation of mass then yields:
\begin{subequations}
\begin{eqnarray}
\frac{\partial\epsilon_s}{\partial t}+\frac{\partial}{\partial x}(u_0 \epsilon_s)&=&0,\label{eq:massconsa}\\
\frac{\partial\epsilon_w}{\partial t}+\frac{\partial}{\partial x}(u_0 \epsilon_w)&=&-\dot m\label{eq:massconsb}.
\end{eqnarray}%
\label{eq:masscons}%
\end{subequations}%
Here, $\dot m$ is the drying rate, with dimensions of $\mathrm{kg}/(\mathrm{m}^3\cdot \mathrm{s})$, and is a general function of $\epsilon_s,\epsilon_l$, and $T$.  
The density of the product phase is denoted by $\overline{\rho}=\epsilon_s+\epsilon_l$, this satisfies the equation
\begin{equation}
\frac{\partial\overline{\rho}}{\partial t}+\frac{\partial}{\partial x}(u_0 \overline{\rho})=-\dot m.
\end{equation}

The energy balance equation is formulated in terms of $e$, the internal energy per unit mass in the product stream.  This is given by:
\begin{equation}
e=\left[\sum_{i=s,l} c_{p,i}(\epsilon_i/\overline{\rho})\right]\left(T-\Tref\right).
\end{equation}
Here, $T$ is the temperature of the product, $\Tref$ is a reference temperature, and the ratios $\epsilon_i/\overline{\rho}$ play the role of mass fractions.  Also, $c_{p,i}$ is the specific heat at constant pressure (per unit mass basis) in the $i^\text{th}$ component of the product mixture. Thus, the equation for energy balance in the product stream reads:
\begin{equation}
\frac{\partial \overline{\rho} e}{\partial t}+\frac{\partial}{\partial x}(u_0 \overline{\rho}e)=\dot q-\dot m h_l.
\label{eq:energycons}
\end{equation}
Here, $\dot q$ is the heat transfer into the product (units: $\mathrm{J}/\left(\mathrm{m^3}\cdot\mathrm{s}\right)$), and $h_l$ is the latent  heat of evaporation, per unit mass.   Equation~\eqref{eq:energycons} can be further simplified to read:
\begin{equation}
\frac{\partial T}{\partial t}+u_0\frac{\partial T}{\partial x}=\frac{\dot q-\dot m \left[h_l-c_{p,w}(T-\Tref)\right]}{\sum_{i=s,l}c_{p,i}\epsilon_{i}},
\label{eq:energycons1}
\end{equation}
which is an advection equation for temperature.
Thus, there are three state variables: the densities $\epsilon_w$ and $\epsilon_l$, and the product temperature $T$.
Two laws  / closure models are required (one for $\dot q$ and one for $\dot m$).   The three  model equations~\eqref{eq:massconsa},~\eqref{eq:massconsb}, and~\eqref{eq:energycons1} transfer inputs at the inlet to outputs at the outlet.  

\subsection{Drying rate}

The functional form of the drying rate has two extreme cases.  In the initial drying stage, the drying rate is limited by the ability of the moist gas to take up evaporated water from the product.  This is the so-called \textit{constant drying-rate regime}.  In this regime, the drying rate reads\cite{mujumdar2006handbook}:
\begin{equation}
\dot m=k_c a_V\left(Y_*(T)-Y_{air}\right).
\end{equation}
 Here, $ k_c$ is a constant with units of $\mathrm{kg}/(\mathrm{m}^2\cdot\mathrm{s})$, $a_V$ is the surface area of the product in contact with the moist gas, per unit volume, $Y_*$ is the absolute humidity of the air at saturation, at the product temperature $T$.  Also, $Y_{air}$ is the absolute humidity of the air.  Here, the absolute humidities are expressed on a dry-gas basis, that is, $Y=M_{wv}/M_{da}$, where $M_{wv}$ is the mass of water vapour in a parcel of moist gas, and $M_{da}$ is the amount of dry gas in that same parcel.  

%The saturation value of absolute humidity $Y_*$ can be computed by going through various standard expressions in thermodynamics, starting with:
%%
%\begin{equation}
%Y_*(T)=0.62198\left(\frac{p_{ws}}{p_a-p_{ws}}\right),
%\end{equation}
%where $p_{ws}$ is the saturation pressure of the water and $p_a$ is the pressure of the moist gas (assumed constant).  The saturation pressure of the water can then be obtained from Antoine's equation,
%%
%\begin{equation}
%p_{ws}=10^{A-B/(C+T)},
%\end{equation}  
%where $A$, $B$, and $C$ are the (known) Antoine's coefficients.  Thus, there is only a single unknown model constant $k_c$ that needs to be applied in the constant drying-rate regime.

Beyond a certain critical moisture level in the solid, there is no longer sufficient water at the product surface to enable direct evaporation.  The critical moisture level is given in terms of the  absolute moisture content of the product on a dry-solid basis, $X=\epsilon_{w}/\epsilon_s$, the critical level is denoted by $X_\mathrm{c}$.
 Below this level, the drying is diffusion-dominated, as it relies on diffusive processes to bring water from the interior of the product to the surface, where evaporation can then occur.  This is the falling drying-rate regime, in this regime the drying rate is given by the Lewis equation\cite{mujumdar2006handbook},
\begin{equation}
\dot m=k_f \epsilon_s\left(X-X_*\right).
\label{eq:lewis1}
\end{equation}
Here, $k_f$ is a rate coefficient (units: $\mathrm{s}^{-1}$), and $X_*$ is the equilibrium moisture content of the solid.  In some applications, $k_f$ is known to be a function of the product temperature.  Equation~\eqref{eq:lewis1} can be rewritten as:
\begin{equation}
\dot m=k_f \left(\epsilon_l-X_*\epsilon_w\right);
\label{eq:lewis2}
\end{equation}
Summarizing, the expression for the drying rate which accounts for the two different regimes is:
\begin{equation}
\dot m=\begin{cases}k_c a_V\left(Y_*(T)-Y_{air}\right),&X>X_\mathrm{c},\\
                    k_f \left(\epsilon_l-X_*\epsilon_w\right),&X\leq X_{\mathrm{c}}.\end{cases}
\label{eq:dotmgen}
\end{equation}
where $k_c$ and $k_f$ are the two rate coefficients which are not known \textit{a priori}.

In the present article, we are concerned mostly with the optimization of drying, which requires functional differentiation of the drying rate.  Hence, the discontinuous drying rate in Equation~\eqref{eq:dotmgen} is not suitable.  In practice, the drying-rate curve (which encodes all information about $\dot m$) exhibits a smooth transition from the constant-rate regime to the falling-rate regime.  Therefore, a more appropriate functional form of the drying rate would be to have a smooth interpolation between the two regimes in Equation~\eqref{eq:dotmgen}.  However, for the purposes of showcasing the optimization techniques which we develop in this article, we instead focus  exclusively on the falling drying-rate regime, and use the expression $\dot m=k_f(\epsilon_l-X_*\epsilon_w)$, with the assumption that both $k_f$ and $X_*$ are constants.

\subsection{Heating}

We consider a constant amount of thermal energy delivered to the product along the length of the drier, this can be achieved for instance by indirect contact of the product with a closed-loop flow of steam.  The thermal energy per unit time, per unit volume is denoted by $\dot q$.  Thus, the total power is
\begin{equation}
P=A_\times\int_0^\length \dot q\,\mathd x,
\end{equation}
where $A_\times$ is the cross-sectional area of the drier.  It is assumed that the power is delivered to the product evenly along the length of the drier, thus, we can write 
\begin{equation}
\dot q=\frac{P}{A_\times \length}.
\end{equation}

%$\dot q=cx$, where $c$ is a constant, hence $c=2P/(A_\times \length^2)$, and:
%%
%\begin{equation}
%\dot q=\frac{2P}{A_\times \length}(x/\length).
%\end{equation}

\section{Industrial Drier: Equilibrium Solutions, Linearization}
\label{sec:linear}

In this section we consider again the three-equation model of the industrial drier, which we summarize here as:
\begin{subequations}
\begin{eqnarray}
\frac{\partial\epsilon_s}{\partial t}+\frac{\partial}{\partial x}(u_0 \epsilon_s)&=&0,\\
\frac{\partial\epsilon_w}{\partial t}+\frac{\partial}{\partial x}(u_0 \epsilon_w)&=&-\dot m,\\
\frac{\partial T}{\partial t}+u_0\frac{\partial T}{\partial x}&=&\underbrace{\frac{\dot q-\dot m \left[h_l-c_{p,w}(T-\Tref)\right]}{\sum_{i=s,l}c_{p,i}\epsilon_{i}}}_{=H},
\end{eqnarray}%
\label{eq:threeeqn}%
\end{subequations}%
where $\dot m=k_f(\epsilon_l-X_*\epsilon_s)$.
We focus first of all on developing an equilibrium solution of this model, and then develop a linear stability analysis, outlining the conditions under which the equilibrium solution is stable to small perturbations at the inlet.  

%Throughout this section, we assume that $k_f$ and $X_*$ are constants, independent of temperature; we discuss the implications of lifting this assumption at the end of the section.

\subsection{Equilibrium Solution}

We consider Equation~\eqref{eq:threeeqn} with fixed inlet conditions $\epsilon_s(0,t)=\epsilon_{s0}$, $\epsilon_w(0,t)=\epsilon_{w0}$, and $T(0,t)=T_0$, where $\epsilon_{s0}$, $\epsilon_{w0}$, and $T_0$ are constants.  Under these inlet conditions, we seek an equilibrium solution of Equation~\eqref{eq:threeeqn}. We immediately obtain $\epsilon_s(x)=\epsilon_{s0}$, the constant solution.  The other two equations in the set simplify:
\begin{subequations}
\begin{eqnarray}
u_0\frac{\mathd \epsilon_w}{\mathd  x}&=&-k_f\left(\epsilon_w-X_*\epsilon_s\right),\\
u_0\frac{\mathd  T}{\mathd  x}&=&H,\label{eq:threeeqn_eqmb}
\end{eqnarray}%
\label{eq:threeeqn_eqm}%
\end{subequations}%
We continue with the assumption that $\dot q$ is constant, and equal to $\dot q=P/(A_{\times} \length)$, where $P$ is the power.
%
%We note briefly that a semi-implicit expression for $\epsilon_w(x)$ can be obtained using the integrating-factor technique:
%
%\begin{equation}
%\epsilon_w(x)=\epsilon_{w0}\mathe^{-\int_0^x [k_f(T)/u_0]\mathd x}+\frac{\epsilon_{s0}}{u_0}\mathe^{-\int_0^x [k_f(T)/u_0]\mathd x}\int_0^x \mathe^{\int_0^{x'} [k_f(T)/u_0]\mathd x''}X_* k_f(T)\mathd x'.
%\end{equation}
%

\begin{figure}[htb]
	\centering
		\includegraphics[width=0.6\textwidth]{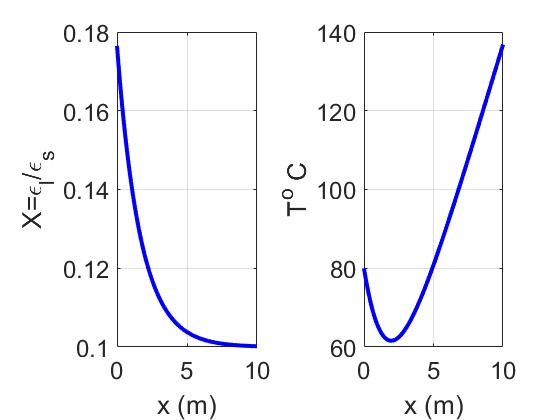}
		\caption{Plot showing the behavior of the equilibrium solution.  The parameter values used to produce the simulation are in Table~\ref{tab:params2}.}
	\label{fig:eqm_example}
\end{figure}

An explicit expression for $\epsilon_w(x)$ can be obtained using the integrating-factor technique:
%This is a fully-explicit expression in the case when $k_f$ and $X_*$ are independent of $T$:
%
\begin{equation}
\epsilon_w(x)=\epsilon_{w0}\mathe^{-k_f x/u_0}+\epsilon_{s0}X_*\left(1-\mathe^{-k_f x/u_0}\right).
\label{eq:validation}
\end{equation}
Then, back-substitution into Equation~\eqref{eq:threeeqn_eqmb} reduces the three-equation model down to a single equation for temperature.
However, for our purposes, it is more straightforward to solve the three-equation model numerically.   Hence, numerical solutions to the three-equation model are shown in Figure~\ref{fig:eqm_example}.
%
%
%where the integrating factor $\mathe^{-\int_0^x [k_f(T)/u_0]\mathd x}$
%
Here, we have used the constant value $k_f=1/(5\, \mathrm{min})$, based on a typical drying time for solids in a dryer.  We further use the typical parameter values in 
Table~\ref{tab:params2} 
%Table~2 
 (the velocity in Table 2 corresponds to a residence time of $30\,\mathrm{mins}$, $T_{R}=\length/u_0$).
\begin{table}[htb]
		\caption{Parameter values for the generation of the equilibrium solution.  The auxiliary parameters $\dot\mu_{inlet}$ and $\Phi_{inlet}$ have been introduced, these then fix $\epsilon_{i,inlet}$.}
		{\begin{tabular}{|c|c|}
		\hline
		Parameter & Value\\
		\hline
		  $P$ & $4\times 10^4\,\mathrm{Watt}$\\
			$\length$ & $10\,\mathrm{m}$\\
			$u_0$ & $1/3\,\mathrm{m}\cdot\mathrm{min}^{-1}$\\
			$X_*$ & 0.1\\
			$\Tinlet$ & $80^\circ\mathrm{C}$\\
			$x_{w,inlet}$ & 0.15\\
		%\hline
		%\end{tabular}
		%
		%
		%\begin{tabular}{|c|c|}
		%\hline
		%Parameter & Value\\
		%\hline
			$x_{s,inlet}$ & $1-x_{w,inlet}$\\
			$\dot\mu_{inlet}$ & $10\,\mathrm{kg}\cdot\mathrm{min}^{-1}$\\
			$A_{\times}$  & $\pi(0.5\,\mathrm{m})^2$\\
			$\Phi_{inlet}$ & $\dot\mu_{inlet}/A_{\times}$\\
			$\overline{\rho}_{inlet}$& $\Phi_{inlet}/u_0$\\
			$\epsilon_{i,inlet}$ & $x_{i,inlet}\overline{\rho}_{inlet}$\\
			$c_{p,s}$ & $1980.4\,\mathrm{J}/(\mathrm{kg}\cdot\mathrm{K})$\\
			$c_{p,l}$ & $4181.5\,\mathrm{J}/(\mathrm{kg}\cdot\mathrm{K})$\\
			$h_l$     & $2.25\times 10^6\,\mathrm{J}/\mathrm{kg}$\\
		\hline
		\end{tabular}}
		\label{tab:params2}
\end{table}

\subsection{P\'eclet Number}

Based on Reference~\cite{zhang2012modelling}, we take a value of $k=0.6\,\mathrm{W}/(\mathrm{kg}\cdot \mathrm{K})$, the same for each component in the product (liquid and solid).  The thermal diffusivity is therefore $k/\sum_{i=s,l}c_{p,i} \epsilon_i$, and the P\'eclet number is:
\begin{equation}
\mathrm{Pe}=\frac{u_0 \ell}{k/\sum_{i=s,l}c_{p,i}\epsilon_i}.
\end{equation}
Based on the values in 
%Table~2
 Table~\ref{tab:params2}, 
and using the values at the inlet for $\epsilon_i$, the P\'eclet number is $8172$.  The large value of the P\'eclet number justifies the neglect of diffusion on the three-equation model.  Furthermore, the absence of sharp temperature gradients in the numerical solutions suggests that there are no `boundary-layer effects', i.e. small regions where the diffusion term is important.  These results justify the choice of a hyperbolic model in the present work.

\subsection{Linear Stability Analysis}

We consider a time-dependent inlet condition
\begin{equation}
\epsilon_{s}(0,t)=\epsilon_{s,0}+\delta \epsilon_{s}(t),\qquad
\epsilon_{l}(0,t)=\epsilon_{l,0}+\delta \epsilon_{l}(t),\qquad
T(0,t)=T_0+\delta T(t),
\label{eq:lsa1}
\end{equation}
where $\delta\epsilon_s(t)$, $\delta\epsilon_{l}(t)$, and $\delta T$ are small-amplitude time-dependent fluctuations.  Hence, the solution of the three-equation model~\eqref{eq:threeeqn} can be decomposed as:
\begin{equation}
\epsilon_s(x,t)=\epsilon_{s,0}+\delta\epsilon_{s}(x,t),\qquad
\epsilon_l(x,t)=\epsilon_{l,eq}(x)+\delta\epsilon_{s}(x,t),\qquad
T(x,t)=T_{eq}(x)+\delta T(x,t),\qquad
\label{eq:lsa2}
\end{equation}
where $\epsilon_{l,eq}(x)$ and $T_{eq}(x)$ denote the equilibrium solution.  We substitute Equation~\eqref{eq:lsa2} into the three-equation model~\eqref{eq:threeeqn} and linearize around the equilibrium solution to obtain:
\begin{equation}
\frac{\partial}{\partial t}\left(\begin{array}{c}\delta\epsilon_s \\ \delta\epsilon_l \\ \delta T\end{array}\right)
+u_0\frac{\partial}{\partial x}\left(\begin{array}{c}\delta\epsilon_s \\ \delta\epsilon_l \\ \delta T\end{array}\right)
=\mathbb{J}\left(\begin{array}{c}\delta\epsilon_s \\ \delta\epsilon_l \\ \delta T\end{array}\right),
\label{eq:lsa3}
\end{equation}
where $\mathbb{J}$ is the Jacobian matrix
\begin{equation}
\mathbb{J}=\left(\begin{array}{ccc}
0 & 0 & 0\\
-\frac{\partial \dot m}{\partial \epsilon_s} & -\frac{\partial\dot m}{\partial \epsilon_l}  & -\frac{\partial \dot m}{\partial T} \\
\frac{\partial H}{\partial \epsilon_s} & \frac{\partial H}{\partial \epsilon_l}  & \frac{\partial H}{\partial T} \\
\end{array}\right)
\label{eq:jac}
\end{equation}
We decompose the inlet perturbations into the harmonic components, using
\begin{equation}
\delta T(t)=\frac{1}{2\pi}\int_{-\infty}^\infty \widehat{\delta T}_{\omega}\mathe^{\imag \omega t}\mathd\omega,
\end{equation}
and similarly for the inlet densities.  The perturbations in the body of the drier (i.e. $x\in(0,\length]$) can also be assumed to have the same decomposition.  Thus, the linear problem~\eqref{eq:lsa3} can be re-written in Fourier space as:
\begin{equation}
u_0\frac{\partial}{\partial x}\left(\begin{array}{c}\widehat{\delta\epsilon}_{s,\omega} \\ \widehat{\delta\epsilon}_{l,\omega} \\ \widehat{\delta T}_\omega\end{array}\right)
=\left(-\imag\omega\mathbb{I}+\mathbb{J}\right)\left(\begin{array}{c}\widehat{\delta\epsilon}_{s,\omega} \\ \widehat{\delta\epsilon}_{l,\omega} \\ \widehat{\delta T}_\omega\end{array}\right)
\label{eq:lsa4}
\end{equation}
The solution of Equation~\eqref{eq:lsa3} can be formally written as:
\begin{equation}
\left(\begin{array}{c}\widehat{\delta\epsilon}_{s,\omega}(x) \\ \widehat{\delta\epsilon}_{l,\omega}(x) \\ \widehat{\delta T}_\omega(x)\end{array}\right)=
\mathe^{[-\imag \omega x + \int_0^x \mathbb{J}\mathd x]/u_0}\left(\begin{array}{c}\widehat{\delta\epsilon}_{s,\omega}(0) \\ \widehat{\delta\epsilon}_{l,\omega}(0) \\ \widehat{\delta T}_\omega(0)\end{array}\right).
\end{equation}

Hence, the question as to whether initial small-amplitude disturbances at the inlet decay to zero as they move downstream is determined by the eigenvalues of $\int_0^x\mathbb{J}\mathd x$.
As this matrix has a row of zeros, one of the eigenvalues is itself zero, therefore the characteristic equation simplifies:
\begin{equation}
\left|\begin{array}{cc}-\int_0^x \frac{\partial \dot m}{\partial \epsilon_l}\mathd x-\lambda & -\int_0^x \frac{\partial \dot m}{\partial T}\mathd x\\
\int_0^x\frac{\partial H}{\partial \epsilon_l} \mathd x& \int_0^x\frac{\partial H}{\partial T}\mathd x-\lambda\end{array}\right|=0.
\label{eq:charpoly}
\end{equation}
We first look at the case where $k_f$ and $X_*$ are independent of $T$, in this case there is a single negative eigenvalue and a single positive eigenvalue; the positive eigenvalue is given by:
%
%\begin{equation}
%\lambda=-\int_0^x \frac{\partial \dot m}{\partial \epsilon_l}\mathd x=-k_f x\leq 0,
%\end{equation}
%and
%
\begin{equation}
\lambda=\int_0^x\frac{\partial H}{\partial T}\mathd x=\int_0^x \left(\frac{\dot m c_{p,w}}{c_{pw}\epsilon_w+c_{p,l}\epsilon_l}\right)\mathd x\geq 0.
\label{eq:lambdaplus}
\end{equation}
%Thus, there is one stable eigenvalue and one unstable eigenvalue.  
Hence,  small-amplitude inlet disturbances are amplified as they are convected downstream to the outlet.  
However, as the drier has finite length, the degree to which the inlet disturbances are amplified is also limited by the growth factor $\mathe^{\max_{[0,\length]}{\lambda/u_0}}$
%, which more explicitly reads as follows:
%
%\begin{equation}
%\text{Max. Amplification}=\exp\left[\frac{k_f}{u_0}\int_0^\length \left(\frac{\epsilon_l-X_*\epsilon_s}{c_{p,l}\epsilon_l+c_{p,s}\epsilon_s}\right)\mathd x\right].
%\label{eq:maxamp}
%\end{equation}
This result applies to all inlet forcing frequencies $\omega$.  Thus, the instability is reminiscent of a convective instability\cite{huerre1985absolute} -- effectively, an amplifier, which amplifies all inlet disturbances at the same rate as they are convected downstream.  This is in contrast to an absolute instability, which is present in systems with a natural frequency.  It also contrasts to the simplified mathematical model in Section~\ref{sec:simple}, which is a stable system, in the sense that disturbances at the inlet are damped to zero as they propagate to the outlet.

As the perturbations away from the inlet are bounded, it is possible to introduce a proportional (bounded) perturbation in the control variable $\dot q$ to bring the outlet temperature $T(\length,t)$ to a desired value.  To show this, we consider again the linear problem~\eqref{eq:lsa3}, however, we now consider further the possibility that the heat source $\dot q$ is also subject to a perturbation, which the operator can control.  Hence, we consider $\dot q=\dot q_0+\delta \dot q$, where $\dot q_0$ is a constant and equal to $(P/A_\times \length)$, and $\delta\dot q$ depends only on time.
%
%$\dot q_0=(2P/A_x\length)(x/\length)$, and $\delta\dot q$ depends on time alone.  Thus, Equation~\eqref{eq:lsa3} becomes
%
\begin{equation}
\frac{\partial}{\partial t}\left(\begin{array}{c}\delta\epsilon_s \\ \delta\epsilon_l \\ \delta T\end{array}\right)
+u_0\frac{\partial}{\partial x}\left(\begin{array}{c}\delta\epsilon_s \\ \delta\epsilon_l \\ \delta T\end{array}\right)
=\mathbb{J}\left(\begin{array}{c}\delta\epsilon_s \\ \delta\epsilon_l \\ \delta T\end{array}\right)+\delta\dot q(t) \left(\begin{array}{c}0\\ 0 \\ \eta(x)\end{array}\right),
\label{eq:lsa_withq}
\end{equation}
where $\eta(x)=\left[c_{p,l}\epsilon_{l}(x)+c_{p,s}\epsilon_{s,0}\right]^{-1}$.  In Fourier space, this reads:
\begin{equation}
u_0\frac{\partial}{\partial x}\left(\begin{array}{c}\widehat{\delta\epsilon}_{s,\omega} \\ \widehat{\delta\epsilon}_{l,\omega} \\ \widehat{\delta T}_\omega\end{array}\right)
=\left(-\imag\omega\mathbb{I}+\mathbb{J}\right)\left(\begin{array}{c}\widehat{\delta\epsilon}_{s,\omega} \\ \widehat{\delta\epsilon}_{l,\omega} \\ \widehat{\delta T}_\omega\end{array}\right)
+\widehat{\delta \dot q}_\omega \left(\begin{array}{c}0\\ 0 \\ \eta(x)\end{array}\right).
\label{eq:lsa_withq_fourierx}
\end{equation}

We look again at Equation~\eqref{eq:lsa_withq_fourierx}, with the same restrictions as before that $k_f$ and $X_*$ are constants.  Then, there is an explicit solution for $\widehat{\delta\epsilon}_{s,\omega}(x)$ and $\widehat{\delta\epsilon}_{l,\omega}(x)$:
\begin{eqnarray*}
\widehat{\delta\epsilon}_{s,\omega}(x)&=&\widehat{\delta\epsilon}_{s,\omega}(0)\mathe^{-\imag\omega x/u_0},\\
\widehat{\delta\epsilon}_{l,\omega}(x)&=&\widehat{\delta\epsilon}_{l,\omega}(0)\mathe^{-\imag\omega x/u_0-k_fx/u_0}+X_*\widehat{\delta\epsilon}_{s,\omega}(0)\mathe^{-\imag\omega x/u_0}\left(1-\mathe^{-k_fx/u_0}\right).
\end{eqnarray*}
Thus, inlet oscillations are propagated undamped to the outlet, this is due to 
%(i) the zero eigenvalue of the Jacobian $\mathbb{J}$ and (ii), the Gramian being singular.
 the zero eigenvalue of the Jacobian $\mathbb{J}$.
Correspondingly, the equation for the temperature fluctuations reads:
\begin{equation}
\frac{\mathd}{\mathd x}\widehat{\delta T}_\omega+\left(\frac{\imag\omega - \frac{\partial H}{\partial T}}{u_0}\right)\widehat{\delta T}_\omega=
\frac{1}{u_0}\underbrace{\left(\frac{\partial H}{\partial \epsilon_s}\widehat{\delta\epsilon}_{s,\omega}+\frac{\partial H}{\partial \epsilon_l}\widehat{\delta\epsilon}_{s,\omega}\right)}_{=\rho(x)}+\widehat{\delta \dot q}_\omega\eta(x)/u_0,
\end{equation}
with solution
\begin{multline}
\widehat{\delta T}_\omega(x)=\widehat{\delta T}_\omega(0)\mathe^{-\imag \omega x/u_0+\lambda_+(x)/u_0}\\
+\mathe^{-\imag \omega x/u_0+\lambda_+(x)/u_0}\int_0^x [\rho(x')/u_0]\mathe^{\imag\omega x'/u_0-\lambda_+(x')/u_0}\mathd x'\\
+\widehat{\delta \dot q}_\omega \mathe^{-\imag \omega x/u_0+\lambda_+(x)/u_0}\int_0^\length [\eta(x')/u_0]\mathe^{\imag\omega x'/u_0-\lambda_+(x')/u_0}\mathd x'.
\end{multline}
Here, $\lambda_+(x)$ is the positive eigenvalue in Equation~\eqref{eq:lambdaplus}.
Hence, in order to suppress the temperature fluctuations at the outlet, such that $\widehat{\delta T}_\omega(\length)=0$, it suffices to apply the control:
%we lo
\begin{equation}
\widehat{\delta \dot q}_\omega=\frac{-\widehat{\delta T}_\omega(0)-\int_0^\length [\rho(x')/u_0]\mathe^{\imag\omega x'/u_0-\lambda_+(x')/u_0}\mathd x'}{\int_0^\length [\eta(x')/u_0]\mathe^{\imag\omega x'/u_0-\lambda_+(x')/u_0}\mathd x'},
\label{eq:oc_analytic}
\end{equation}
which exists, provided the denominator is non-zero.  

To account for the hyperbolic nature of the system, great care must be taken in reconstructing the solution of $\delta T(x,t)$ in the spatio-temporal domain, from the data given in $\widehat{\delta T}_\omega(x)$.  In particular, the reconstruction needs to take account of the jump which occurs at $x=u_0 t$, where the solution $\delta T(x,t)$ transitions from being dominated by information at the initial time, to information from the inlet.   This was straightforward in the case of the simplified model (Section~\ref{sec:simple}), but is considerably more difficult here because the forcing term $\delta \dot q(t)\eta(x)$ depends on space and time.  We do not pursue this approach any further here, and instead pass over to a numerical approach.

\subsection{Numerical Results}

\begin{figure}
	\centering
		\includegraphics[width=0.45\textwidth]{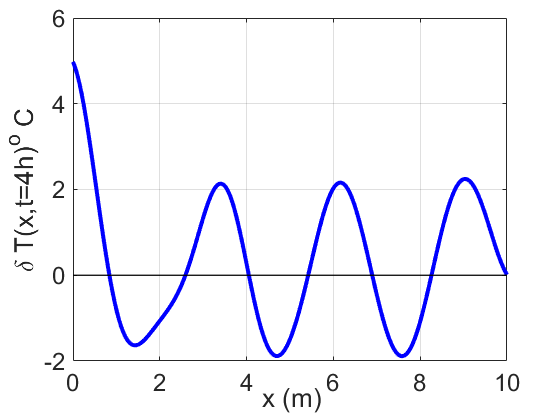}
		\includegraphics[width=0.45\textwidth]{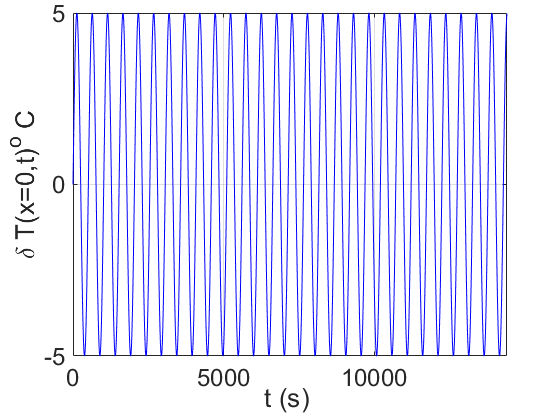}
		\includegraphics[width=0.45\textwidth]{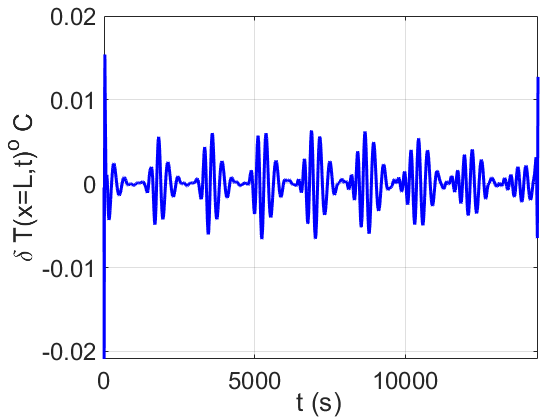}
		\includegraphics[width=0.45\textwidth]{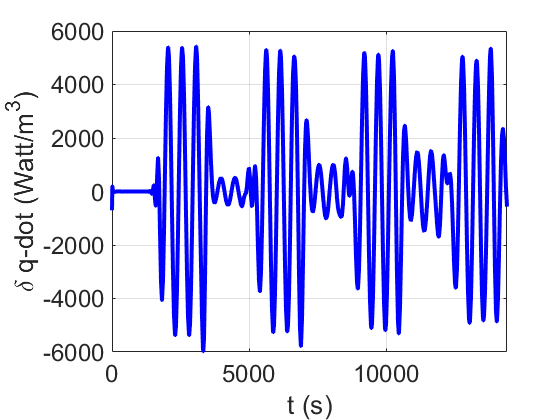}
		\caption{Optimal control theory for the linear system.  (a) $\delta T(x,t)$ at the simulation final time 
		$t=4\,\mathrm{hr}$; (b) $\delta T(x,t)$ at the inlet $x=0$, for all times; (c) the controlled temperature fluctuation $\delta T(\length,t)$ showing attainment of the control $\delta T(\length,t)\rightarrow 0$; (d) the optimal control.}
	\label{fig:lsa_control}
\end{figure}

We apply optimal control theory to the linear system~\eqref{eq:lsa3}, with a view to keeping the outlet temperature fluctuation at zero, $\delta T(\length,t)=0$.  We use the numerical technique developed in Section~\ref{sec:simple}, for which the corresponding adjoint problem reads:
\begin{equation}
-\frac{\partial}{\partial t}\left(\begin{array}{c}\psi_s \\ \psi_l \\ \psi_T\end{array}\right)
-u_0\frac{\partial}{\partial x}\left(\begin{array}{c}\psi_s \\ \psi_l \\ \psi_T\end{array}\right)
-\mathbb{J}^T\left(\begin{array}{c}\psi_s \\ \psi_l \\ \psi_T\end{array}\right)=0,
\end{equation}
with terminal condition $\psi(x,\tau)=0$ at the end-time $\tau$ and outlet conditions
\begin{equation}
\psi_s(\length,t)=0,\psi_l(\length,t)=0,\qquad \psi_T(\length,t)=-(1/u_0)\delta T(\length,t).
\end{equation}
The $n^{\mathrm{th}}$ guess for the optimal control is labelled by $[\delta\dot q(t)]^n$, and the updated guess at the next iteration is given by
\begin{equation}
[\delta\dot q(t)]^{n+1}=[\delta\dot q(t)]^n+\lambda\int_0^\length \Psi_T(x,t) \eta(x,t)\mathd x,
\end{equation}
where $\eta(x,t)=\left[c_{p,l}\epsilon_{l}(x,t)+c_{p,s}\epsilon_{s}(x,t)\right]^{-1}$, and
where $\lambda$ is a small parameter that is chosen using the Barzilai--Borwein method.
We look at  monochromatic inlet forcing
\begin{equation}
\delta T(t)=(5^\circ\mathrm{C})\sin(\omega_F t),
\label{eq:deltaT_forcing}
\end{equation}
where $\omega_F=2\pi/T_F$, and the forcing period $T_F$ is taken to be $8.5\,\mathrm{min}$.  Also, the inlet solid and liquid density fluctuations are set to zero.
Numerical results are shown in Figure~\ref{fig:lsa_control} for $1,000$ iterations of the steepest-descent algorithm, with a final residual value $0.0378$, this can be reduced further by applying further iterations of the steepest-descent algorithm.
We  use $N=200$ gridpoints and $\Delta t=0.1\,\mathrm{s}$, for converged numerical results.   The optimization is carried out over a four-hour time horizon (i.e. $\tau=4\,\mathrm{hr}$).

\begin{figure}
	\centering
		\includegraphics[width=0.6\textwidth]{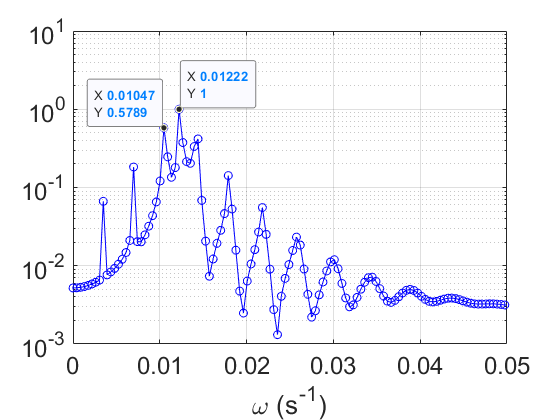}
		\caption{Normalized power spectrum of $\delta \dot q$ taken with respect to a simulation over 4 hours.  The maximum at $\omega\approx \omega_F$ is highlighted, as well as a secondary maximum.}
	\label{fig:lsa_control_fourierx}
\end{figure}

We focus our attention on the behaviour of $\delta\dot q$, shown in Figure~\ref{fig:lsa_control}(d).  This consists of a single dominant oscillation with frequency $\omega_F$, i.e. to match the inlet forcing.  There is a second, long-time trend present, in the form of a `beat', as evidenced by the long-time periodic modulation of the optimal control $\delta\dot q$ in the figure.  To understand this in more detail, we have computed the power spectrum of the optimal control, shown in Figure~\ref{fig:lsa_control_fourierx}.
There is a single sharp peak in the data at $\omega\approx \omega_F$ (specifically, $\omega=0.01222$) and a secondary peak nearby, this time at $\omega=0.01047$.  
The secondary peak is due to the `ramp-up' in the forcing term which occurs at $x=u_0 t$, where the solution goes from initial-condition-dominated to inlet-dominated.  This is confirmed by looking at other (non-optimal) forcing terms $\delta \dot q$ -- only those with the `ramp-up' exhibit a secondary peak.
As these two peaks are close in frequency, this gives rise to a `beating' phenomenon: the beat frequency is $\omega_B=0.01222-0.01047$ (no division by two), which corresponds to a beat period of $1\,\mathrm{hours}$, this manifests itself in the real-space signal as the long-period modulation.  

\section{Industrial Drier: Nonlinear Optimal Control Theory}
\label{sec:nonlinear}

We now apply optimal control theory to the nonlinear system~\eqref{eq:threeeqn}, with a view to keeping the outlet temperature fluctuation at a set point $T(\length,t)=T_*$.  We use the numerical technique developed in Section~\ref{sec:simple}.  Hence, we solve the nonlinear problem forward in time with a guess $\dot q(t)$ for the optimal forcing, specifically, we solve:
\begin{eqnarray*}
\frac{\partial\epsilon_s}{\partial t}+\frac{\partial}{\partial x}(u_0 \epsilon_s)&=&0,\\
\frac{\partial\epsilon_w}{\partial t}+\frac{\partial}{\partial x}(u_0 \epsilon_w)&=&-\dot m,\\
\frac{\partial T}{\partial t}+u_0\frac{\partial T}{\partial x}&=&\underbrace{\frac{\dot q-\dot m \left[h_l-c_{p,w}(T-\Tref)\right]}{\sum_{i=s,l}c_{p,i}\epsilon_{i}}}_{=H}.
\end{eqnarray*}
Next, we solve the adjoint problem:
\begin{equation}
-\frac{\partial}{\partial t}\left(\begin{array}{c}\psi_s \\ \psi_l \\ \psi_T\end{array}\right)
-u_0\frac{\partial}{\partial x}\left(\begin{array}{c}\psi_s \\ \psi_l \\ \psi_T\end{array}\right)
-\mathbb{J}^T\left(\begin{array}{c}\psi_s \\ \psi_l \\ \psi_T\end{array}\right)=0,
\end{equation}
with terminal condition $\psi(x,\tau)=0$ at the end-time $\tau$ and outlet conditions
\begin{equation}
\psi_s(\length,t)=0,\psi_l(\length,t)=0,\qquad \psi_T(\length,t)=-(1/u_0)[ T(\length,t)-T_*].
\end{equation}
Here, the Jacobian $\mathbb{J}$ has exactly the same form as in the linear stability analysis (Section~\ref{sec:linear}), Hence,  although the forward problem is nonlinear, the adjoint problem is linear.
However,  the Jacobian must now be evaluated on the densities and temperature from the forward-problem, which vary in space and time.  Consequently, knowledge of the forward variables at all points in space and time is required for the adjoint computation -- this means that the simulation time is much longer than before in the case of the fully linear problem.   These are  generic features of optimal control problems for partial differential equations.     As before, we label the $n^{\mathrm{th}}$ guess for the optimal control  by $[\dot q(t)]^n$, and the updated guess at the next iteration is given by
\begin{equation}
[\dot q(t)]^{n+1}=[\dot q(t)]^n+\lambda\int_0^\length \psi_T(x,t) \eta(x)\mathd x,
\end{equation}
where $\lambda$ is a small parameter that is chosen (again) using the Barzilai--Borwein method.  This approach is iterated until the outlet temperature is sufficiently close the set point, for all times.

\subsection{Results}

For the initial condition of the forward problem, we use the equilibrium solution from Section~\ref{sec:linear}, $\epsilon_{s,0}$ (a constant), $\epsilon_{l,eq}(x)$, and $T_{eq}(x)$.  The inlet conditions are then selected as:
\begin{subequations}
\begin{eqnarray}
\epsilon_s(0,t)=\epsilon_{s,0},\\
\epsilon_l(0,t)=\epsilon_{l,eq}(0)\left[1+\delta\alpha \sin(\omega_F t)\right],\\
T(0,t)=T_{eq}(0)\left[1+\delta\alpha \sin(\omega_F t)\right],
\end{eqnarray}%
\label{eq:ics}%
\end{subequations}%
where $\omega_F$ is the forcing frequency, which we again take to be $\omega_F=2\pi/(8.5\,\mathrm{mins})$.   We take the magnitude of the initial inlet perturbation $\delta\alpha$ to be $0.05$.  Also, in order to avoid a discontinuous control at $t=0$, we also take the set point $T_*$ to be $T_{eq}(\length)$. 

\begin{figure}
	\centering
		\includegraphics[width=0.45\textwidth]{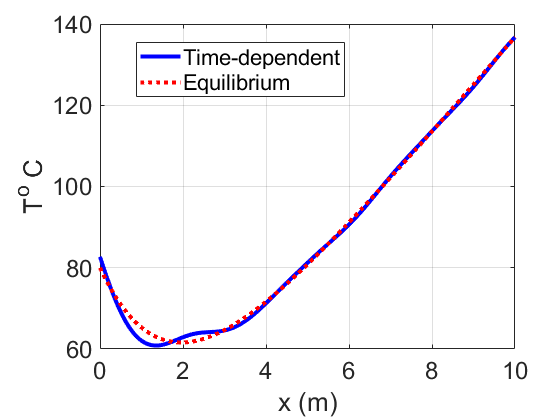}
		\includegraphics[width=0.45\textwidth]{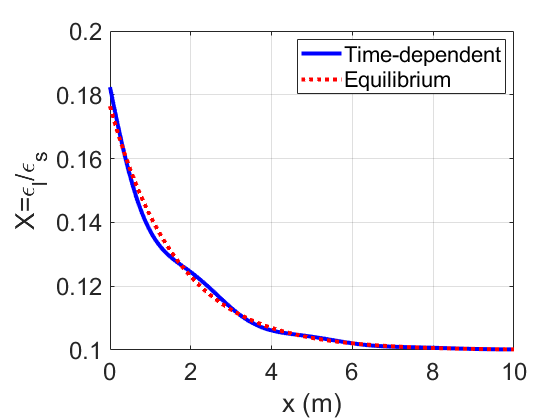}
		\includegraphics[width=0.45\textwidth]{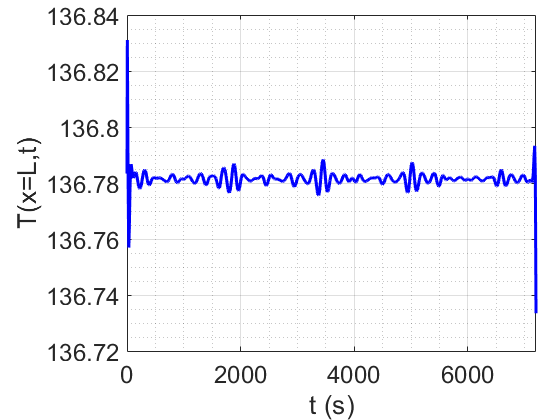}
		\includegraphics[width=0.47\textwidth]{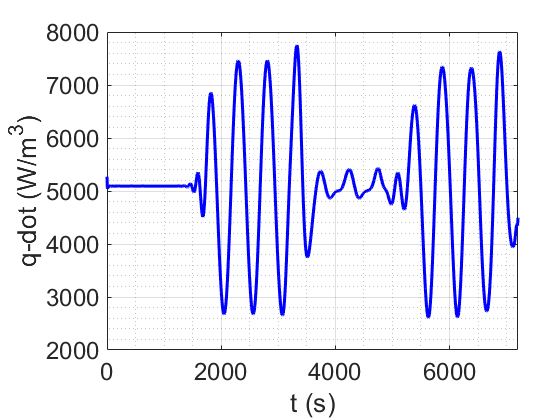}
		\caption{Optimal control theory for the nonlinear system, $\delta\alpha=0.05$.  (a) $T(x,t)$ at the simulation final time 
		$t=2\,\mathrm{hr}$; (b) Moisture content $X(x,t)$ at the simulation final time(c) the controlled temperature fluctuation $T(\length,t)$ showing attainment of the control $\delta T(\length,t)\rightarrow T_*$; (d) the optimal control.}
	\label{fig:nonlin1}
\end{figure}

\begin{figure}
	\centering
		\includegraphics[width=0.6\textwidth]{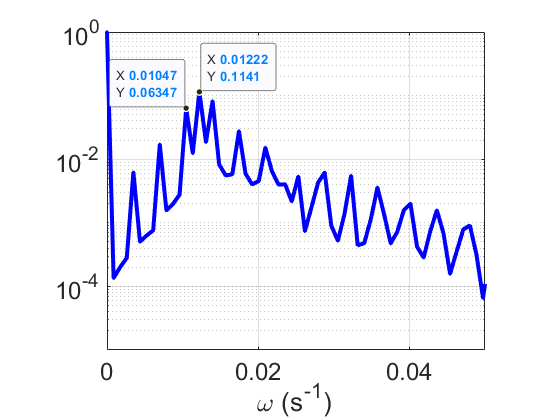}
		\caption{Normalized power spectrum of $\delta \dot q$ taken with respect to a simulation over 2 hours (nonlinear case).  The maximum at $\omega\approx \omega_F$ is highlighted.}
	\label{fig:nonlin2}
\end{figure}

Results are shown in Figure~\ref{fig:nonlin1}, for a simulation  out to $\tau=2\,\mathrm{h}$.
These are qualitatively very similar to the previous linear results.  The optimal control $\dot q(t)$ keeps the outlet temperature at the set point $T_*$.  A power-spectrum analysis of $\dot q(t)$ reveals a single peak at the forcing frequency, as well as neighboring peaks, thus producing a `beat' pattern in the optimal forcing (Figure~\ref{fig:nonlin2}) -- qualitatively the same as in the linear case.  The large peak in the power spectrum in Figure~\ref{fig:nonlin2} at $\omega=0$ corresponds to the constant component, which is treated separately in the linear stability analysis.

\begin{figure}
	\centering
		\includegraphics[width=0.45\textwidth]{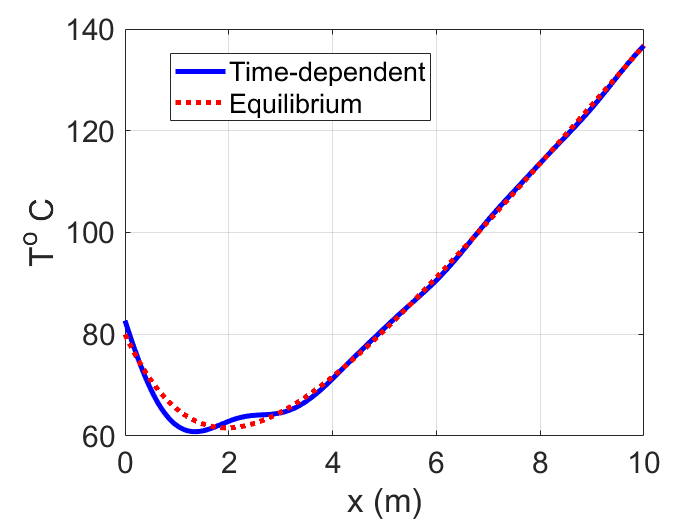}
		\includegraphics[width=0.45\textwidth]{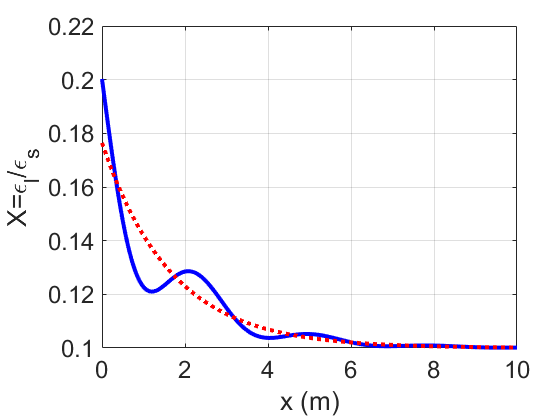}
		\includegraphics[width=0.45\textwidth]{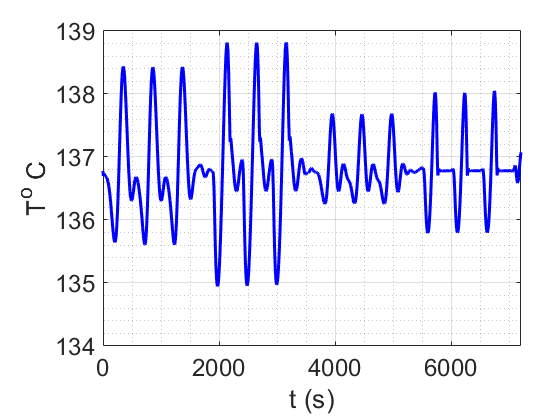}
		\includegraphics[width=0.44\textwidth]{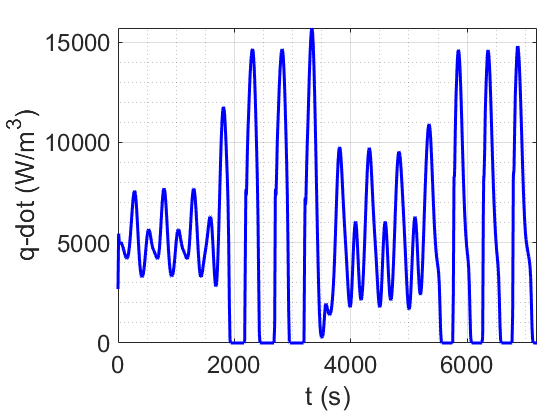}
		\caption{Optimal control theory for the nonlinear system, large-amplitude case with $\delta\alpha=0.4$  (a) $T(x,t)$ at the simulation final time 
		$t=2\,\mathrm{hr}$; (b) Moisture content $X(x,t)$ at the simulation final time; (c) the outlet temperature $T(\length,t)$ showing the deviation of $T(\length,t)$ from the set-point $T_*$, which is minimized under (d), the optimal control.}
	\label{fig:nonlin1_lg}
\end{figure}

We have also explored large-amplitude initial disturbances using the same initial condition as before (Equation~\eqref{eq:ics}), but with $\delta\alpha=0.2$.  In this case, the optimal control has $q(t)\leq 0$ for some $t$, which corresponds to cooling the product, as opposed to heating.  As this is technically infeasible, we look at imposing a constraint $q(t)\geq 0$ for all $t\in[0,\tau]$, this can be achieved for instance by letting $q(t)=(1/2)[\theta(t)]^2$, and performing the steepest-descent calculations using $\delta \lagrange/\delta\theta=\theta(t) \int_0^\length \psi_T(x,t)\eta(x)\mathd x$.
Results are shown in Figure~\ref{fig:nonlin1_lg}, and are now qualitatively different from the previous linear and small-amplitude nonlinear cases.  As we are now constrained to have $\dot q\geq 0$, the set point $T(\length,t)=T_*$ is not reached, instead, the difference between $T(\length,t)$ and $T_*$ is minimized in the root-mean-square sense.
However, the optimal control in Figure~\ref{fig:nonlin1_lg}(c) still  exhibits the previous `beat' pattern.  
%%A power spectrum analysis (Figure~\ref{fig:nonlin2_lg})
%%
%%
%\begin{figure}
	%\centering
		%\includegraphics[width=0.6\textwidth]{ps_nonlin_lg}
		%\caption{Normalized power spectrum of $\delta \dot q$ taken with respect to a simulation over 2 hours (nonlinear case, large amplitude).  The local maximum at $\omega\approx \omega_F$ is highlighted, as well as the second-harmonic generation at $2\omega_F$ -- this is now the dominant peak in the spectrum.}
	%\label{fig:nonlin2_lg}
%\end{figure}
%reveals that the dominant mode in the optimal control (apart from the constant contribution) is now in fact $2\omega_F$, and that $3\omega_F$ is also important, i.e. second- and third-harmonic generation.  This is a straightforward consequence of nonlinear interaction between the various modes in the problem. 

%In the general case, the exchange of energy between the normal modes (such as in Figure~\ref{fig:nonlin2_lg}) may lead to chaos and hence, to late-time uncertainty in the results.  Consequently,  caution is warranted when using the optimal control theory in the nonlinear regime for long periods of time, if the aim is to control an underlying physical system.  In this scenario, the basic optimal control methodology may need to be corrected, using model-predictive control.  %However, in the present case, the spectral peak at $2\omega_F$ is so well 

\subsection{Outlet Moisture Level}

From Figures~\ref{fig:nonlin1}--\ref{fig:nonlin1_lg} (in particular, panel (b) in both figures), it can be seen that the outlet moisture level $X(\length,t)$ remains close to the equilibrium moisture content $X_*$.  This is a simple consequence of our choice of parameters for the drier -- the velocity $u_0$, and the length of the dryer have both been chosen with reference to the value of $k_f$ such the outlet equilibrium moisture level $X_{eq}(\length,t)$ is close to $X_*$.  Furthermore, inspection of the linear stability analysis in Section~\ref{sec:linear} shows that the liquid's density fluctuations are damped by the exponential factor $\mathe^{-k_fx/u_0}$ where again, the value of $\length$ means that significant damping of fluctuations occurs along the length of the channel.  Thus, the effective control of the outlet moisture levels can be thought of as byproduct of good design.  
%This is discussed in more detail in Section~\ref{subsec:discussion}.   
In contrast, the attainment of the control $T(\length,t)\rightarrow T_*$ is a consequence of active control whereby $\dot q(t)$ is selected according to the steepest-descent method.

\section{Discussion and Conclusions}
\label{sec:conc}

Summarizing, we have developed a mathematical framework for spatio-temporal control in industrial drying, starting from first principles, namely convection-type models for heat and mass transfer.    In contrast to previous models, ours take account of the temperature variations of product which can occur inside the drier itself, these can be significant.
We have applied optimal control theory to the models, using a Steepest-Descent Method with the Barzilai--Borewein formula for stepsize, to determine the optimal external heating required to maintain the outlet temperature close to a set-point-temperature.  For the three-equation model, and in the case of small-amplitude inlet forcing, the optimal control maintains the outlet temperature at the set point; for large amplitudes this is not possible, and the optimal control minimizes the root-mean-square difference between the outlet temperature and the set-point temperature.

The main aim of the present work is to control the outlet temperature of the product, using a time-dependent heat source  as a single control variable. The outlet moisture level has been maintained at the fixed, equilibrium value by virtue of a judicious choice of design parameters -- including the residence time and the drier length.  In cases where these parameters cannot be chosen in this way, a second manipulated variable should be considered -- for instance, a valve at the outlet, allowing for the product to stay in the drier for longer.  In this case, the general mathematical framework introduced in the present work still applies.

The optimal control theory requires advance knowledge of the inlet conditions. 
In practical settings, this may not be available; however, knowledge of the state of the plant obtained from the model will be supplemented by sensor data for temperature at various points along the drier.  In this case, a Kalman Filter can be used to correct the model, in order to keep the model predictions `close' to observations.  The Kalman Filter may include process noise -- which may be appropriate for modelling inlet conditions that consist of a mean value plus fluctuations.  This approach can be supplemented with model predictive control -- a key  methodology from control theory used for estimating the optimal control based on the given mathematical model, the sensor data coming from the plant, and the Kalman Filter.
Our work provides a theoretical basis for understanding how model predictive control would work in case of an industrial drier with temporal and spatial variations.  First, the model predictive control relies on a robust discretized numerical model -- such as the ones developed herein.  Second, the optimal controls presented herein 
 still form a `best case scenario' against which the performance of an online model predictive control can be assessed.  Hence, the models and methods developed herein form a key validation step which may be developed further in a more practical setting.

\subsection*{Acknowledgements}

The author acknowledges insights from discussions with       James Mc Carthy of Cervos Ltd., as well as from Luis Carnevale and Miguel Bustamante.  The author finally acknowledges the work done by Odhr\'an Dooley in helping to compile the parameter values in Table~\ref{tab:params2}.

\subsection*{Funding}

This work has been produced as part of ongoing work within the ThermaSMART network.
The ThermaSMART network has received funding from the European Union's Horizon 2020
research and innovation programme under the Marie Sklodowska--Curie grant agreement No.
778104.

%\bibliographystyle{unsrt}
%\bibliography{drying_bibliography}

\end{document}